\documentclass{amsart}
\usepackage{amssymb}
\usepackage{amsfonts, amsmath, amsthm,graphics,amssymb}
\newsavebox{\uuunit}
\sbox{\uuunit}
    {\setlength{\unitlength}{0.825em}
     \begin{picture}(0.6,0.7)
        \thinlines
        \put(0,0){\line(1,0){0.5}}
        \put(0.15,0){\line(0,1){0.7}}
        \put(0.35,0){\line(0,1){0.8}}
       \multiput(0.3,0.8)(-0.04,-0.02){12}{\rule{0.5pt}{0.5pt}}
     \end {picture}}
\newcommand{\ID}{\mathord{\!\usebox{\uuunit}}}

\def\tr{{\rm Tr}}
\def\nn{\nonumber}

\theoremstyle{plain}
\newtheorem{theorem}{Theorem}[section]
\newtheorem{exam}[theorem]{Example}
\newtheorem{lm}[theorem]{Lemma}
\newtheorem{thm}[theorem]{Theorem}
\newtheorem{cor}[theorem]{Corollary}

\theoremstyle{definition}
\newtheorem{df}[theorem]{Definition}
\newtheorem{rmk}[theorem]{Remark}

\begin{document}

\title[Reductions of the Schlesinger equations]{On the Reductions and Classical Solutions of the Schlesinger equations.}

\author{B. Dubrovin}
\address{SISSA, International School of Advanced
Studies, via Beirut 2-4, 34014 Trieste, Italy}
\email{dubrovin@sissa.it}
\author{M. Mazzocco}
\address{School of Mathematics, The University of Manchester, Manchester M60 1QD, United Kingdom.}
\email{Marta.Mazzocco@manchester.ac.uk}
\thanks{
The authors wish to dedicate this paper to the memory of A. Bolibruch who was a great source of inspiration.  We are grateful to B. Malgrange, H. Umemura and M. Singer for many helpful conversations. The researches of B.D. were
partially supported by Italian Ministry of Education research grant Cofin2004
``Geometry of Integrable Systems'', and also by the ESF research grant MISGAM.  The researches of M.M. was supported by EPSRC,
SISSA, IRMA (Strasbourg) and the European network ENIGMA. Finally, the authors wish to thank the referee for carefully reading the manuscript and for her/his helpful comments.}

\maketitle

\hfill{\it To the memory of our friend Andrei Bolibruch}
\vskip 2cm

\begin{abstract} The Schlesinger equations  $S_{(n,m)}$
describe monodromy preserving deformations of order $m$ Fuchsian systems
with $n+1$ poles. They can be considered as a family of commuting time-dependent
Hamiltonian systems on the direct product of $n$ copies of $m\times m$
matrix algebras equipped with the standard linear Poisson bracket.
In this paper we address the problem of reduction of particular solutions 
of ``more complicated'' Schlesinger equations $S_{(n,m)}$ to ``simpler''   $S_{(n',m')}$
having $n'< n$ or $m' < m$. 
\end{abstract}

\tableofcontents

\section{Introduction.} 

The {\it Schlesinger equations} $S_{(n,m)}$ \cite{Sch} is the following 
system of nonlinear differential equations 
\begin{eqnarray}
&&
{\partial\over\partial u_j} {A}_i= 
{[ {A}_i, {A}_j]\over u_i-u_j},\qquad i\neq j,\nonumber \\
&&
{\partial\over\partial u_i} {A}_i= 
-\sum_{j\neq i}{[ {A}_i, {A}_j]\over u_i-u_j},\label{sch}
\end{eqnarray} 
for $m\times m$ matrix valued functions ${A}_1=A_1(u),\dots,
{A}_{n}=A_n(u)$, where the independent variables $u=(u_1, \dots, u_n)$ must be 
pairwise distinct. The first non-trivial case $S_{(3,2)}$
of the Schlesinger equations corresponds to the 
famous sixth Painlev\'e equation \cite{fuchs, Sch, Gar1}, 
the most general of all Painlev\'e
equations. In the case of any number $n>3$ of $2\times2$ 
matrices ${A}_j$, the Schlesinger equations reduce to the
Garnier systems ${\mathcal G}_n$ (see \cite{Gar1,Gar2,Ok1}).

The Schlesinger equations $S_{(n,m)}$ appeared in the theory of 
{\it isomonodromic deformations}\/
of Fuchsian systems. Namely, the monodromy matrices of the Fuchsian system
\begin{equation}
{d \Phi\over  d{z}} =\sum_{k=1}^{n}{{A}_k(u)\over {z}-u_k}\Phi,
\qquad\qquad {z}\in{\mathbb C}\backslash\{u_1,\dots,u_{n}\}
\label{N1}\end{equation}
do not depend on $u=(u_1, \dots, u_n)$ if the matrices $A_i(u)$ satisfy
(\ref{sch}). Conversely, under certain assumptions on the 
matrices $A_1$, \dots, $A_n$ and for the matrix
\begin{equation}\label{a-inf}
A_\infty:= -\left( A_1+\dots +A_n\right),
\end{equation}
{\it all} isomonodromic deformations of the Fuchsian system are given by
solutions to the Schlesinger equations (see, e.g.,
\cite{sib})\footnote{Bolibruch constructed non-Schlesinger isomonodromic 
deformations in \cite{Bol}. These can occur when the matrices $A_i$ are resonant,
i.e. admit pairs of eigenvalues with positive integer differences. To avoid such non-Schlesinger isomonodromic deformations, we need to extend the set of monodromy data (see Section 2 below).}.

The solutions to the Schlesinger equations can be parameterized by the
{\it monodromy data} of the Fuchsian system (\ref{N1}) (see precise definition below 
in Section 2). To reconstruct the solution starting from given monodromy 
data one is to solve the classical {\it Riemann - Hilbert problem} of 
reconstruction of the Fuchsian system
from its monodromy data. The main outcome of this approach says that
the solutions $ {A}_i(u)$ can be continued 
analytically to meromorphic functions on the universal covering of
$$
\left\{(u_1,\dots,u_{n})\in{\mathbb C}^{n}\,|\,u_i\neq u_j\,\hbox{for}\, 
i\neq j\right\}
$$
\cite{Mal,Mi}.
This is a generalization of the celebrated {\it Painlev\'e
property}\/ of absence of movable critical singularities (see details in
\cite{ince,its}). In certain cases the technique based on the theory of
Riemann--Hilbert problem gives a possibility to compute the asymptotic 
behavior
of the solutions to the Schlesinger equations near the critical locus $u_i=u_j$
for some $i\neq j$, although, in general, the problem of determining the asymptotic behaviuor 
near the critical points is still open
\cite{jimbo, DM, guzzetti, costin}. 

It is the Painlev\'e property that was used by Painlev\'e and Gambier
as the basis for their classification scheme of nonlinear differential
equations. Of the list of some 50 second order nonlinear differential equations
possessing Painlev\'e property
the six (nowadays known as {\it Painlev\'e equations}) 
are selected due to the following crucial property: the general solutions
to these six equations cannot be expressed in terms of {\it classical
functions}, i.e., elementary functions, elliptic and other classical
transcendental functions (see \cite{Um1} for a modern approach to this
theory based on a nonlinear version of the differential Galois theory). 
In particular, according to these results
the general solution to the Schlesinger system $S_{(3,2)}$ corresponding 
to Painlev\'e-VI equation cannot be expressed in terms of classical functions.

A closely related question is the problem of construction and classification
of {\it classical solutions} to Painlev\'e equations and their generalizations.
This problem remains open even for the case of Painlev\'e-VI although there are
interesting results based on the theory of symmetries of Painlev\'e equations
\cite{OK, OK3, kitaev} and on the geometric approach to studying the space 
of monodromy
data \cite{DM, Hit, M,M2}.

The above methods do not give any clue to solution of 
the following general problems: are 
solutions of $S_{(n+1,m)}$ or of $S_{(n,m+1)}$ more complicated than those
of $S_{(n,m)}$? Which solutions to $S_{(n+1,m)}$ or $S_{(n,m+1)}$
can be expressed via solutions to $S_{(n,m)}$? Furthermore, which of them can
ultimately be expressed via classical functions?

In this paper we aim to suggest a general approach to the theory of reductions and classical solutions 
of the general Schlesinger equations $S_{(n,m)}$ for all $n$, $m$, based on the Riemann--Hilbert 
problem and on the group-theoretic properties of the monodromy group of the linear system 
(\ref{N1}). Our approach consists in determining the monodromy data of the Fuchsian system 
(\ref{N1}) that guarantee to have a reduction to $S_{(n-1,m)}$ or $S_{(n,m-1)}$ and eventually a
classical solution. 

We need a few definitions.  Let us fix a solution to the Schlesinger equations $S_{(n,m)}$.
Applying the algebraic operations and differentiations to the matrix
entries of this solution 
we obtain a field ${\mathcal S}_{(n,m)}$
equipped with $n$ pairwise commuting differentiations $\partial/\partial u_1$, \dots,
$\partial / \partial u_n$, to be short a {\it differential field}. Define the {\it rational
closure} ${\mathcal K}$ of a differential field ${\mathcal S}$ represented by functions of $n$
variables by taking rational functions with coefficients in ${\mathcal S}$
$$
{\mathcal K}:= {\mathcal S}(u_1, \dots, u_n).
$$
Taking the rational closure of the differential field $S_{(n,m)}$, we obtain the differential 
field ${\mathcal K}_{(n,m)}$. (Needless to say that the field ${\mathcal K}_{(n,m)}$ depends on the choice of the solution to the Schlesinger equations $S_{(n,m)}$.)

We now construct new differential fields obtained
from ${\mathcal K}_{(n_1, m_1)}$, \dots, ${\mathcal K}_{(n_k, m_k)}$ by applying one or
more
of the following {\it admissible} elementary operations.

1. Tensor product. Given two differential fields ${\mathcal K}_1$ and ${\mathcal K}_2$
represented by functions of $n_1$ and $n_2$ variables $u_1$, \dots, $u_{n_1}$
and $v_1$, \dots, $v_{n_2}$ respectively, produce
new differential field ${\mathcal K}_1 \otimes{\mathcal K}_2$ taking the rational
closure of the minimal differential field of functions of $n_1+n_2$ 
independent variables $u_1$, \dots, $u_{n_1}$, $v_1$, \dots, $v_{n_2}$ containing both $K_1$ and $K_2$. A particular case of 
this operation is

2. Addition of an independent variable. Given a differential field ${\mathcal K}$
represented by functions of $n$ variables $u_1$, \dots, $u_n$ define an
extension
$\tilde{\mathcal K}\supset {\mathcal K}$ by taking rational functions of a new independent
variable $u_{n+1}$ with coefficients in ${\mathcal K}$ ,
$$
\tilde{\mathcal K}= {\mathcal K}\otimes {\mathbb C}(u_{n+1}).
$$

3. Given two differential fields ${\mathcal K}_1$, ${\mathcal K}_2$ represented by
functions of the same number of variables $n$, define the {\it composite} 
${\mathcal K}_1 {\mathcal K}_2$ taking the minimal differential field of 
functions of $n$ variables containing both
${\mathcal K}_1$ and ${\mathcal K}_2$ and applying the rational closure procedure.

4. A differential field extension ${\mathcal K}'\supset {\mathcal K}$ is said to be of
the {\it Picard--Vessiot type} if it is the minimal rationally closed
differential field of functions of $n$
variables containing solutions of a Pfaffian linear system with coefficients
in ${\mathcal K}$ \cite{kol,singer}. 

Recall that a Pfaffian linear system of the order $k$ with coefficients
in a differential field ${\mathcal K}$ represented by functions of $n$ variables
reads
$$
{\partial Y\over \partial u_i} = B_i Y, \quad i=1, \dots, n, \quad Y=(y_1, \dots, y_k)^T
$$
where the $n$ matrices
$$
B_i \in Mat(k; {\mathcal K}) \quad {\rm for }\quad i=1,
\dots, n
$$
must satisfy
$$
{\partial B_i\over \partial u_j} -{\partial B_j\over \partial u_i} +[B_i, B_j]=0, \quad {\rm
for}~ {\rm all}~ i\neq j.
$$
The linear space of solutions of the Pfaffian system is finite dimensional. The  differential field ${\mathcal K}'$ is the minimal extension of ${\mathcal K}$ containing all
components $y_1$, \dots, $y_k$ of any of these solutions.
 
We will also denote ${\mathcal K}^{(N)}$ the differential field extension of 
${\mathcal K}$ obtained by $N$ Picard--Vessiot type extensions of
${\mathcal K}$
$$
{\mathcal K} \subset {\mathcal K}' \subset {\mathcal K}'' \subset 
\dots \subset {\mathcal K}^{(N)}.
$$

Using the above admissible extensions we can describe in what circumstances
a particular solution to the Schlesinger equations $S_{(n,m)}$ can 
be expressed via
solutions to $S_{(n', m')}$ with smaller $n'$ or $m'$. Similar results were
obtained in \cite{M1} for the special case of $m=2$.

\begin{thm} \label{itro1}Consider a solution to $S_{(n,m)}$ such that 
the eigenvalues of the matrix $A_\infty$ are pairwise distinct and
the monodromy group of the associated Fuchsian system (\ref{N1}) admits a
$k$-dimensional invariant subspace, $k>0$. Then this solution belongs
to a Picard--Vessiot type extension ${\mathcal K}^{(N)}$ for some $N$
of the composite 
$$
{\mathcal K} = {\mathcal K}_{(n,k)} {\mathcal K}_{(n,m-k)}
$$
 where ${\mathcal K}_{(n,k)}$ and ${\mathcal K}_{(n,m-k)}$ are two differential fields associated with certain two solutions of the systems
${\mathcal S}_{(n,k)}$ and ${\mathcal S}_{(n,m-k)}$ respectively.
\end{thm}
In particular,
\begin{cor}\label{tria}
Given a solution to $S_{(n,m)}$ such that the monodromy group of the 
associated Fuchsian system (\ref{N1}) is upper-triangular  and the eigenvalues of $A_\infty$ are 
pairwise distinct, it
belongs to a Picard--Vessiot type extension 
${\mathcal K}_0^{(N)}$ for some $N$ of
$$
{\mathcal K}_0={\mathbb C}(u_1,\dots,u_n).
$$
\end{cor}

The proof of this Theorem is based on the following two lemmata.

\begin{lm}\label{perv} Given a solution
\begin{equation}\label{sist-a}
A(z; u) =\sum_{i=1}^n \frac{A_i(u)}{z-u_i}
\end{equation}
to the Schlesinger equations ${\mathcal S}_{(n,m)}$ with diagonalizable matrix $A_\infty$ such that the associated monodromy representation has a $k$-dimensional invariant subspace, denote ${\mathcal K}_{(n,m)}$ the corresponding differential field. Then there exists a matrix
$$
G(z;u) \in \bar{\mathcal K}_{(n,m)}[z], \quad \det G(z) \equiv 1
$$
such that all matrices $B_i(u)$, $i=1, \dots, n$ of the gauge equivalent Fuchsian system with
\begin{equation}\label{sist-b}
B(z;u) =G^{-1}(z;u) A(z;u) G(z;u) + G^{-1}(z;u) \frac{dG(z;u)}{dz}=\sum_{i=1}^n \frac{B_i(u)}{z-u_i}
\end{equation}
have a $u$-independent  $k$-dimensional common invariant subspace. Here $\bar{\mathcal K}_{(n,m)}$ is a Picard - Vessiot type extension of the field ${\mathcal K}_{(n,m)}$. Moreover,
the matrices $B_1(u)$, \dots, $B_n(u)$ satisfy Schlesinger equations.
\end{lm}

This lemma, apart from polynomiality of the gauge transformation in $z$, is the main result of the papers \cite{malek1,malek3} by S.Malek\footnote{Actually, there is a stronger claim in the main result of \cite{malek3}, namely, it is said that the coefficients of the reducing gauge transformation are rational functions in $u_1$, \dots, $u_n$ and entries of $A_1(u)$, \dots, $A_n(u)$. We were unable to reproduce this result.}. We give here a new short proof of this result (for the sake of technical simplicity we add the assumption of diagonalizability of the matrix $A_\infty$) by presenting a
reduction algorithm consisting of  a number of elementary and explicitly written transformations.

It is a one-line calculation that shows that the ${\mathcal S}_{(n,m)}$ Schlesinger equations for the matrices $B_1(u)$, \dots, $B_n(u)$ of the form
$$
B_i(u) =\left(\begin{array}{cc} B_i'(u) & C_i(u)\\ 0 & B_i''(u)\end{array}\right),
$$
where $B_i'(u)$ and $B_i''(u)$ are respectively $k\times k$ and $(m-k)\times (m-k)$ matrices, reduces to the ${\mathcal S}_{(n,k)}$
and ${\mathcal S}_{(n,m-k)}$ Schlesinger systems for the matrices
$B_i'(u)$ and $B_i''(u)$ and to the linear Pfaffian equations
\begin{eqnarray}
&&
\partial_j C_i = \frac1{u_i-u_j} \left( B_i' C_j - B_j' C_i +C_i B_j'' -C_j B_i''\right), \quad j\neq i
\nonumber\\
&&
\partial_i C_i = -\sum_{j\ne i}\frac1{u_i-u_j} \left( B_i' C_j - B_j' C_i +C_i B_j'' -C_j B_i''\right).
\nonumber
\end{eqnarray}
Therefore the Schlesinger deformation of the reduced system (\ref{sist-b}) belongs to a Picard - Vessiot type extension of the composite ${\mathcal K}_{n,k}{\mathcal K}_{n,m-k}$.

To complete the proof of Theorem \ref{itro1} we need to invert the above gauge transformation, i.e., to express the coefficients of the original Fuchsian system (\ref{sist-a}) via the solution of the reduced system (\ref{sist-b}). 

\begin{lm} \label{vtor} (i) For a Fuchsian system (\ref{sist-a}) satisfying the assumptions of the previous lemma,  the monodromy data, in the sense of Definition \ref{daty} here below,
$$
\Lambda^{(1)}(A), R^{(1)}(A), \dots, \Lambda^{(\infty)}(A), R^{(\infty)}(A),  C_1(A), \dots, C_n(A)
$$
of the system (\ref{sist-a}) and
$$
\Lambda^{(1)}(B), R^{(1)}(B), \dots, \Lambda^{(\infty)}(B), R^{(\infty)}(B),  C_1(B), \dots, C_n(B)
$$
of (\ref{sist-b})
are related by
\begin{eqnarray}\label{sviaz}
&&
\Lambda^{(i)}(B)=P^{-1}\Lambda^{(i)}(A)P, \quad i=1, \dots, n
\nn\\
&&
R^{(i)}(B)=P^{-1}R^{(i)}(A)P, \quad i=1, \dots, \infty, 
\\
&&
C^{(i)}(B)=P^{-1}C^{(i)}(A)P, \quad i=1, \dots, n
\nn\\
&&
\nn\Lambda^{(\infty)}(B)=P^{-1}\Lambda^{(\infty)}(A)P+ {\rm diag}\,(N_1, \dots, N_m), \quad N_i\in{\mathbb Z}.
\end{eqnarray}
Here $P\in S_m$ is a permutation matrix.

\noindent (ii) Denote ${\mathcal K}_{n,m}^A$ and ${\mathcal K}_{n,m}^B$ the differential fields associated with the Schlesinger deformations of two systems (\ref{sist-a}) and (\ref{sist-b}) respectively.
If the monodromy data of the systems are related as in (\ref{sviaz}) then there exists a matrix
$$
\tilde G(z; u) \in \bar {\mathcal K}_{n,m}^B[z], \quad \det \tilde G(z; u)\equiv 1
$$
such that
$$
A(z) \equiv \tilde G^{-1}(z;u) B(z;u) \tilde G(z;u) + \tilde G^{-1}(z;u) \frac{d\tilde G(z;u)}{dz}.
$$
Here, like in Lemma \ref{perv},
$\bar {\mathcal K}_{n,m}^B$ is a suitable Picard - Vessiot type extension of the field ${\mathcal K}_{n,m}^B$.
\end{lm}

We obtain therefore two inclusions
\begin{equation}\label{inclu}
{\mathcal K}_{n,m}^B\subset \bar {\mathcal K}_{n,m}^A, \quad {\mathcal K}_{n,m}^A\subset \bar {\mathcal K}_{n,m}^B.
\end{equation}
The Theorem \ref{itro1} easily follows from the above statements.
\medskip

Let us proceed now to the second mechanism of reducing the Schlesinger equations.
Let us assume that $l$ monodromy matrices $M_{i_1}$, \dots, $M_{i_l}$ of the  Fuchsian system of the form (\ref{N1}), are {\it scalar matrices} (i.e., they are proportional to 
$\ID$).
In that case we will call the solution $A_1(u)$, \dots $A_n(u)$  
is {\it $l$-smaller}. We call {\it $l$-erased} the Fuchsian system
${\mathcal S}_{n-l,m}$ of the same size with the poles $z=u_{i_1}$, \dots, $z=u_{i_l}$ erased.

\begin{thm}\label{l-malo}
Let ${A}_{1},\dots,{A}_{n}$ be a $l$--smaller solution of the Schlesinger 
equations.
Then ${A}_{1}(u),\dots , {A}_{n}(u)$ belong to the differential field obtained by 
admissible extensions from ${\mathcal K}_{(n-l, m)}$, the rational closure of 
the differential field ${\mathcal S}_{n-l,m}$ associated with
a solution to the $l$-erased Fuchsian system  $S_{(n-l,m)}$.
In particular, if $l=n-2$ then ${A}_{1},\dots,{A}_{n}$ belong to
the differential field obtained by 
admissible extensions from ${\mathbb C}(u_1,\dots,u_n)$.
\label{small1}\end{thm}

The proof of this Theorem consists first in observing that, due to the fact that all matrices 
$A_1,\dots,A_n$ can be assumed to be traceless, any scalar matrix $M_k$ must have the form
$$
M_k= e^{2\pi i p\over m} \ID, \quad p\in {\mathbb Z}.
$$
As a first step we assume $M_k$, say for $k=n$, to be the identity and we construct a gauge transformation in a suitable Picard - Vessiot type extension of the field ${\mathcal K}_{n,m}^A$ defined by the $l$--solution ${A}_{1},\dots,{A}_{n}$, which maps $A_n$ to zero without changing the nature of the other singular points $u_1,\dots,u_{n-1}$, nor introducing new ones. In this way we obtain a new solution $B_1,\dots,B_{n-1}$ of the Schlesinger equations $S_{n-1,m}$. We then prove that the original solution $A_1,\dots,A_n$ can be constructed in terms of $B_1,\dots,B_{n-1}$ by means of admissible operations.

When $M_k=e^{2\pi i p\over m} \ID$ is not the identity, we need to map  $A_1,\dots,A_n$ bi--rationally to a new solution $\tilde A_1,\dots,\tilde A_n$ of the Schlesinger equations with $\tilde M_k=\ID$. To this end we apply the  birational 
canonical transformations of Schlesinger equations found in \cite{DM1}\footnote{An alternative way, as it was proposed by the referee, would be to replace our canonical transformations by a combination of Schlesinger transformations of \cite{MJ1} with scalar shifts  instead. However,
the birationality of the proposed transformation need to be justified in the resonant case.}.

To present here this class of transformations let us briefly remind
the canonical Hamiltonian formulation of Schlesinger equations $S_{(n,m)}$ of \cite{DM1}.

Recall \cite{MJU,Man1} that Schlesinger equations can be written as Hamiltonian
systems on the Lie algebra
$$
{\mathfrak g}:= \oplus_{i=1}^n gl(m) \ni (A_1, \dots, A_n)
$$
with respect to the standard linear Lie - Poisson bracket on ${\mathfrak g}$
with some quadratic time-dependent Hamiltonians of the form
\begin{equation}\label{hjm}
H_k:=\sum_{l\neq k}{{\rm tr}\left(A_k A_l\right)\over u_k-u_l}.
\end{equation}
Because of isomonodromicity they can be restricted onto the symplectic leaves
$$
{\mathcal O}_1 \times \dots \times {\mathcal O}_n \in {\mathfrak g}
$$
obtained by fixation of the conjugacy classes ${\mathcal O}_1$,\dots,
${\mathcal O}_n$
of the matrices $A_1$, \dots,
$A_n$. The matrix $A_\infty$ given in (\ref{a-inf}) is a common integral of the Schlesinger
equations. Applying the procedure of symplectic reduction \cite{mars-wein}
we obtain the reduced symplectic space 
\begin{eqnarray}\label{level}
&&
\left\{ A_1\in {\mathcal O}_1 , \dots, A_n\in {\mathcal O}_n,\,
A_\infty = \hbox{given diagonal matrix}\right\}
\nonumber\\
&&
\qquad\hbox{modulo simultaneous diagonal conjugations}.
\end{eqnarray}
The dimension of this reduced symplectic leaf in the generic situation
is equal to $2g$ where
$$
g={m(m-1)(n-1)\over2}-(m-1).
$$
In \cite{DM1}  a new system of the so-called {\it isomonodromic Darboux
coordinates} $q_1$, \dots, $q_g$, $p_1$, \dots, $p_g$
on generic symplectic manifolds (\ref{level}) was constructed and the new Hamiltonians were expressed in these coordinates. Let us explain this construction.

The Fuchsian system (\ref{N1}) can be reduced to a scalar differential
equation of the form
\begin{equation}\label{skal}
y^{(m)} = \sum_{l=0}^{m-1} d_l(z) y^{(l)}.
\end{equation}
For example, one can eliminate last $m-1$ components of the vector function
$\Phi$ to obtain a $m$-th order equation for the first component
$y:=\Phi_1$. (Observe that the reduction procedure depends on the choice of
the component of $\Phi$.)
The resulting Fuchsian equation will have regular singularities
at the same points $z=u_1$, \dots, $z=u_n$, $z=\infty$. It will also have
other singularities produced by the reduction procedure. However, they will be
{\it apparent} singularities, i.e., the solutions to (\ref{skal}) will be
analytic in these points. Generically there will be exactly $g$ apparent
singularities (cf. \cite{oht}; a more precise result about the number
of apparent singularities working also in the nongeneric situation was obtained
in \cite{Bol2}); they are the first part $q_1$, \dots, $q_g$ of the
canonical coordinates. The conjugated momenta are defined by
$$
p_i = {\mathrm{Res}}_{z=q_i} \left( d_{m-2}(z) +{1\over 2} d_{m-1}^2(z)\right), \quad i=1,
\dots, g.
$$

\begin{thm}\cite{DM1} Let the eigenvalues of the matrices $A_1$, \dots, $A_n$, $A_\infty$
be pairwise distinct.
Then the map
\begin{equation}\label{theor}
\left\{\begin{matrix} {\rm Fuchsian} ~ {\rm systems} ~ {\rm with} ~{\rm given}
~{\rm poles,}\\
 {\rm given} ~ {\rm eigenvalues} ~ {\rm of} ~ A_1, \dots, A_n,
A_\infty\\
{\rm modulo} ~ {\rm diagonal} ~ {\rm conjugations} \end{matrix}\right\} \to
(q_1, \dots, q_g, p_1, \dots, p_g)
\end{equation}
gives a system of rational
Darboux coordinates on the generic reduced symplectic leaf (\ref{level}).
The Schlesinger equations $S_{(n,m)}$ in these coordinates are written in the
canonical Hamiltonian form
\begin{eqnarray}
&&
{\partial q_i\over \partial u_k} = {\partial{\mathcal H}_k\over \partial p_i}
\nonumber\\
&&
{\partial p_i\over \partial u_k} = -{\partial {\mathcal H}_k\over \partial q_i}\nonumber
\end{eqnarray}
with the Hamiltonians
$$
{\mathcal H}_k = {\mathcal H}_k(q, p; u) = -{\mathrm{Res}}_{z=u_k} \left( d_{m-2}(z) +{1\over 2} d_{m-1}^2(z)
\right), \quad k=1, \dots, n.
$$
\end{thm}

Here {\it rational Darboux coordinates} means that the elementary
symmetric functions $\sigma_1(q)$, \dots, $\sigma_g(q)$ and 
$\sigma_1(p)$, \dots, $\sigma_g(p)$ are rational functions of the coefficients
of the system and of the poles $u_1$, \dots, $u_n$. Moreover, there exists
a section of the map (\ref{theor}) given by rational functions
\begin{equation}\label{sect}
A_i = A_i(q, p), \quad i=1, \dots, n,
\end{equation}
symmetric in $(q_1, p_1)$, \dots, $(q_g, p_g)$ with coefficients depending on $u_1,\dots,u_n$ and on the eigenvalues if the
matrices $A_i$, $i=1, \dots, n, \infty$. All other Fuchsian systems with the same poles $u_1,\dots,u_n$, the same eigenvalues and the same
$(p_1,\dots,p_g,q_1,\dots,q_g)$ are obtained by simultaneous diagonal conjugation
$$
A_i(q,p)\mapsto C^{-1}A_i(q,p)C, \quad i=1, \dots, n, 
\quad C={\rm diag}\, (c_1,
\dots, c_m).
$$

\begin{thm}\cite{DM1}
The Schlesinger equations $S_{(n,m)}$ written in the canonical form of Theorem 
1.3 admit  a group of birational canonical 
transformations $\langle S_2,\dots, S_m,S_\infty\rangle$
\begin{equation}
S_k:\quad \left\{\begin{array}{l} 
\tilde q_i= u_1+u_k-q_i, \quad i=1,\dots,g,\\
\tilde p_i=-p_i,\quad i=1,\dots,g,\\
\tilde u_l = u_1+u_k-u_l,\quad l=1,\dots,n, \\
\tilde {\mathcal H}_l=-{\mathcal H}_l,\quad l=1,\dots,n,
\end{array}\right.
\label{simmk}\end{equation}
\begin{equation}
S_\infty:\quad \left\{\begin{array}{l} 
\tilde q_i={1\over q_i-u_1}, \quad i=1,\dots,g,\\
\tilde p_i= -p_iq_i^2-{2m^2-1\over m}q_i,\quad i=1,\dots,g,\\
\tilde u_l= {1\over u_l-u_1},\quad l=2,\dots,n, \\
u_{1}\mapsto\infty,\\
\infty\mapsto u_{1},\\
\tilde H_1=H_1,\\
\tilde {H}_l= -H_l (u_l-u_1)^2 + (u_l-u_1) (d^0_{m-1}(u_l-u_1))^2 -\\
\quad-(u_l-u_1) {(m-1)(m^2-m-1)\over m}d^0_{m-1}(u_l-u_1),\\
\hbox{for}\quad l=2,\dots,n\\
\end{array}\right.
\label{simminf}\end{equation}
where 
$$
d^0_{m-1}(u_k) = \sum_{s=1}^g {1\over u_k-q_s}  - {m\, (m-1)\over 2} 
\sum_{l\neq k} {1\over u_k-u_l}.
$$
The transformation $S_k$ acts on the monodromy matrices as follows
\begin{eqnarray}
&&
\tilde  M_1=M_1^{-1}.....M_{k-1}^{-1} M_k M_{k-1}....M_1, \nonumber\\
&&
\tilde M_j=M_{j-1},\quad j=2,\dots,k,\nonumber\\
&&
\tilde M_i=M_i,\quad i=k+1,\dots,n.\nonumber
\end{eqnarray}
The transformation $S_\infty$ acts on the monodromy matrices as follows
\begin{eqnarray}
&&
\tilde M_\infty=e^{-{2\pi i \over m}} M_1,\qquad
\tilde M_1= e^{2\pi i \over m}M_\infty,\nonumber\\
&&
\tilde M_j=M_1^{-1}M_j M_1\quad\hbox{for}\quad j=2,\dots,n.\nonumber
\end{eqnarray}
\end{thm}

To conclude, 
Theorems 1.1 and 1.2 show that for certain very special monodromy groups  the 
Schlesinger equations $S_{(n,m)}$ {\it reduce}\/ to solutions of $S_{(n',m')}$ with 
$n'< n$ and/or $m' < m$ . We do not know any other general mechanism of reducibility of Schlesinger equations. As generically the monodromy group of the system
(\ref{N1}) is not reducible nor smaller, we expect that generic solutions of the 
Schlesinger equations $S_{(n,m)}$ 
do not belong to any admissible extension of composites of the differential fields of the form
${\mathcal K}_{(n',m')}$ with $n' < n$ and/or $m' < m$.  The proof of this fact, that is the proof
of {\it irreducibility} of the Schlesinger equations, is still a rather intriguing open problem.



\setcounter{equation}{0}
\setcounter{theorem}{0}

\section{Schlesinger equations as monodromy preserving deformations of Fuchsian
systems.}

In this section we establish our notations, remind a few basic definitions and
prove some technical lemmata that will be useful throughout this paper.

The Schlesinger equations ${\mathcal S}_{(n,m)}$ describe monodromy preserving 
deformations of 
Fuchsian systems (\ref{N1}) with $n+1$ regular singularities at 
$u_1,\dots,u_{n}$, $u_{n+1}=\infty$:
\begin{equation}\label{f-ur}
{{\rm d}\over{\rm d}{z}} \Phi=\sum_{k=1}^{n}{{A}_k\over {z}-u_k}\Phi,
\qquad\qquad {z}\in{\mathbb C}\backslash\{u_1,\dots,u_{n}\}
\end{equation}
${A}_k$ being $m\times m$ matrices independent of ${z}$, 
and $u_k\neq u_l$ for $k\neq l$, $k,l=1,\dots,n+1$. Let us explain 
the precise meaning of this claim.

\subsection{Levelt basis near a logarithmic singularity and local monodromy
data}

A system
\begin{equation} \label{loga}
{d\Phi\over dz} ={A(z)\over z-z_0} \Phi
\end{equation}
is said to have a {\it logarithmic}, or {\it Fuchsian} singularity at $z=z_0$ if
the $m\times m$ matrix valued function $A(z)$ is analytic in some neighborhood of $z=z_0$.
By definition the {\it local monodromy data}\/ of the system is the class of
equivalence of such systems w.r.t. local  gauge transformations
\begin{equation}\label{gauge0}
A(z) \mapsto G^{-1}(z) A(z)\, G(z) + (z-z_0) G^{-1}(z)\partial_z G(z)
\end{equation}
analytic near $z=z_0$ satisfying
$$
\det G(z_0)\neq 0.
$$
The local monodromy can be obtained by choosing a suitable
fundamental matrix solution of the system (\ref{loga}). The most general
construction of such a fundamental matrix was given by Levelt \cite{Levelt}.
We will briefly recall this construction in the form suggested in \cite{Dub1}.

Without loss of generality one can assume that $z_0=0$. Expanding the system
near $z=0$ one obtains
\begin{equation}\label{sys-exp}
{d \Phi\over d z} = \left( {A_0\over z} + A_1 + z\, A_2 + \dots
\right)  \Phi.
\end{equation}
Let us now describe the structure of local monodromy data. 

Two linear operators $\Lambda$, $R$ acting in the complex $m$-dimensional space $V$
$$
\Lambda, \, R : V \to V
$$
are said
to form an {\it admissible pair} if the following conditions are fulfilled.

\noindent 1. The operator $\Lambda$ is semisimple and the operator $R$ is
nilpotent.

\noindent 2. $R$ commutes with $e^{2\pi i \Lambda}$,
\begin{equation}\label{lev1}
e^{2\pi i \Lambda} R = R\, e^{2 \pi i \Lambda}.
\end{equation}
Observe that, due to the last condition the operator $R$ satisfies
\begin{equation}\label{lev2}
R(V_\lambda) \subset \oplus_{k\in {\mathbb Z}} V_{\lambda+k} \quad {\rm for}
\quad {\rm any} \quad \lambda \in {\rm Spec} \, \Lambda,
\end{equation}
where $V_\lambda \subset V$ is the subspace of all eigenvectors of 
$\Lambda$ with
the eigenvalue $\lambda$. The last condition says that

\noindent 3. The sum in the r.h.s. of (\ref{lev2}) contains only 
non-negative values of
$k$.

A decomposition
\begin{equation}\label{lev3}
R=R_0 + R_1 +R_2 + \dots
\end{equation}
is defined where
\begin{equation}\label{lev4}
R_k(V_\lambda) \subset V_{\lambda+k} \quad {\rm for}
\quad {\rm any} \quad \lambda \in {\rm Spec} \, \Lambda.
\end{equation}
Clearly this decomposition contains only finite number of terms.
Observe the useful identity
\begin{equation}\label{lev5}
z^\Lambda R \, z^{-\Lambda} = R_0 + z\, R_1 + z^2 R_2 + \dots.
\end{equation}

\begin{thm}\label{lev} For a system (\ref{sys-exp}) with a logarithmic singularity at
$z=0$ there exists a fundamental matrix solution of the form
\begin{equation}\label{lev6}
\Phi(z) = \Psi(z) z^\Lambda z^R
\end{equation}
where $\Psi(z)$ is a matrix valued function analytic near $z=0$ satisfying
$$
\det \Psi(0)\neq 0
$$
and $\Lambda$, $R$ is an admissible pair.
\end {thm}

The formula (\ref{lev6}) makes sense after fixing a branch of logarithm 
$\log z$ near $z=0$.
Note that $z^R$ is a polynomial in $\log z$ due to nilpotency of $R$.

The proof can be found in \cite{Levelt} (cf. \cite{Dub1}). Clearly $\Lambda$
is the semisimple part of the matrix $A_0$; $R_0$ coincides with its nilpotent
part. The remaining terms of the expansion appear only in the {\it resonant
case}, i.e., if the difference between some eigenvalues of $\Lambda$ is a
positive integer. In the important particular case of a diagonalizable matrix
$A_0$,
$$
T^{-1}A_0 T=\Lambda ={\rm diag}\, (\lambda_1, \dots, \lambda_m)
$$
with some nondegenerate matrix $T$, the matrix function $\Psi(z)$ in the 
fundamental matrix solution (\ref{lev6}) can be obtained in the form
$$
\Psi(z) = T \left( \ID + z \,\Psi_1 + z^2 \Psi_2 + \dots \right).
$$
The matrix coefficients $\Psi_1$, $\Psi_2$, \dots of the expansion as well as the 
components $R_1$, $R_2$, \dots of the matrix $R$ (see (\ref{lev3}))
can be found recursively from the equations
$$
[\Lambda, \Psi_k]-k\, \Psi_k =-B_k+R_k +\sum_{i=1}^{k-1} \Psi_{k-i} R_i - B_i \Psi_{k-i},
\quad k\geq 1.
$$
Here
$$
B_k:= T^{-1} A_k T, \quad k\geq 1.
$$
If $k_{\rm max}$ is the maximal integer among the differences
$\lambda_i-\lambda_j$ then 
$$
R_k =0 \quad {\rm for} ~ k>k_{\rm max}.
$$
Observe that vanishing of the logarithmic terms in the fundamental matrix
solution (\ref{lev6}) is a constraint imposed only on the first $k_{\rm max}$
coefficients $A_1$, \dots, $A_{k_{\rm max}}$ of the expansion (\ref{sys-exp}).

\begin{exam}\label{example-psi} For the Fuchsian system (\ref{N1}) having diagonal the matrix
$$
A_\infty=-(A_1+\dots +A_n)={\rm diag}\, (\lambda_1, \dots, \lambda_m)
$$
the fundamental matrix of Theorem \ref{lev} has the form
$$
\Phi=\left( 1+\frac{\Psi_1}{z} +O\left(\frac1{z^2}\right)\right) z^{-\Lambda} z^{-R},
$$
where
$$
\Lambda=A_\infty, \quad R=R_1 + R_2 +\dots,
$$
\begin{eqnarray}\label{primer}
&&
\left(R_1\right)_{i\, j} =\left\{\begin{array}{cc} \left( B_1\right)_{i\, j}, & \lambda_i=\lambda_j+1\\ 0, & \mbox{\rm otherwise}\end{array}\right.
\nn\\
&&
\nn\\
&&
B_1 =-\sum_k A_k u_k
\nn\\
&&
\left( \Psi_1\right)_{i\, j} =\left\{ \begin{array}{cc} -\frac{\left( B_1\right)_{i\, j}}{\lambda_i-\lambda_j-1}, & 
\lambda_i\neq \lambda_j+1\\  & \\
\mbox{\rm arbitrary}, & \mbox{\rm otherwise}\end{array}\right.
\nn\\
&&
\\
&&
\left(R_2\right)_{i\, j} =\left\{\begin{array}{cc} \left( B_2 -\Psi_1 R_1 +B_1\Psi_1\right)_{i\, j} , & \lambda_i=\lambda_j+2\\ 0, & \mbox{\rm otherwise}\end{array}\right.
\nn\\
&&
\nn\\
&&
B_2 =-\sum_k A_k u_k^2\nn\\
&&
\left( \Psi_2\right)_{i\, j} =\left\{ \begin{array}{cc} \frac{\left( -B_2+\Psi_1 R_1 -B_1\Psi_1\right)_{i\, j}}{\lambda_i-\lambda_j-2}, & 
\lambda_i\neq \lambda_j+2\\  & \\
\mbox{\rm arbitrary}, & \mbox{\rm otherwise}\end{array}\right.
\nn
\end{eqnarray}
etc.
\end{exam}

It is not difficult to describe the ambiguity in the choice of the admissible
pair of matrices
$\Lambda$, $R$ describing the local monodromy data 
of the system (\ref{sys-exp}). Namely, the diagonal matrix $\Lambda$ is defined up to
permutations of diagonal entries. Assuming the order fixed, the ambiguity in the
choice of $R$ can be described as follows \cite{Dub1}.
Denote ${\mathcal C}_0(\Lambda)\subset GL(V)$ the subgroup consisting of invertible linear operators $G:V\to V$ satisfying
\begin{equation}\label{gr-c0}
 z^\Lambda G \,z^{-\Lambda} = G_0+z\, G_1 +z^2 G_2+\dots .
\end{equation}
The definition of the subgroup can be reformulated \cite{Dub1} in terms of invariance of certain flag in $V$ naturally associated with the semisimple operator $\Lambda$.
The matrix $\tilde R$ obtained from $R$ by the conjugation of the form
\begin{equation}\label{eq-r}
\tilde R= G^{-1} R \, G
\end{equation}
will be called {\it equivalent} to $R$. Multiplying (\ref{lev6}) on the right
by $G$ one obtains another  fundamental matrix solution to the same system of the same
structure
$$
\tilde\Phi(z):=\Psi(z)z^\Lambda z^RG =\tilde \Psi(z) z^\Lambda z^{\tilde R}
$$ 
i.e., $\tilde \Psi(z)$ is analytic at $z=0$ with $\det \tilde\Psi(0)\neq 0$.

The columns of the fundamental matrix (\ref{lev6}) form a distinguished basis
in the space of solutions to (\ref{sys-exp}).

\begin{df} The basis given by the columns of the matrix (\ref{lev6}) is called
{\it Levelt basis} in the space of solutions to (\ref{sys-exp}). 
The fundamental
matrix (\ref{lev6}) is called {\it Levelt fundamental matrix solution}.
\end{df}

The monodromy transformation of the Levelt fundamental matrix solution
reads
\begin{equation}\label{lev7}
\Phi\left(z\, e^{2\pi i}\right) = \Phi(z) M, \quad M= 
e^{2\pi i \Lambda} e^{2\pi i R}.
\end{equation}

To conclude this Section let us denote ${\mathcal C}(\Lambda,R)$ the subgroup of invertible transformations of the form
\begin{equation}\label{lev9}
{\mathcal C}(\Lambda,R)=\{ \, G\in GL(V)\, |\, z^\Lambda  G \,  z^{-\Lambda} = \sum_{k\in {\mathbb Z}} G_k z^k \, \mbox{and} \, [G,R]=0\}.
\end{equation}
The subgroups ${\mathcal C}(\Lambda,R)$ and ${\mathcal C}(\Lambda,\tilde R)$
associated with equivalent matrices $R$ and $\tilde R$ are conjugated.
It is easy to see that this subgroup coincides with the centralizer
of the monodromy matrix (\ref{lev7})
\begin{equation}\label{lev8}
G\in {\mathcal C}(\Lambda, R) \quad {\rm iff} \quad G \, 
e^{2\pi i \Lambda} e^{2\pi i R}  = e^{2\pi i \Lambda} e^{2\pi i R} G, \quad \det
G\neq 0.
\end{equation}

Denote 
\begin{equation}\label{lev10}
{\mathcal C}_0(\Lambda, R) \subset {\mathcal C}(\Lambda, R)
\end{equation}
the subgroup consisting of matrices $G$ such that the expansion (\ref{lev9})
contains only non-negative powers of $z$. Multiplying the Levelt fundamental matrix
(\ref{lev6}) by a matrix $G\in {\mathcal C}_0(\Lambda, R)$ one obtains another
Levelt solution to (\ref{sys-exp})
\begin{equation}\label{lev11}
\Psi(z) z^{\Lambda} z^R G = \tilde \Psi(z) z^{\Lambda} z^R.
\end{equation}

In the next Section we will see that the quotient ${\mathcal C}(\Lambda, R)/
{\mathcal C}_0(\Lambda, R)$ plays an important role in the theory of monodromy
preserving deformations.

\subsection{Monodromy data and isomonodromic deformations of a Fuchsian system}

Denote $\lambda^{(k)}_j$, $j=1,\dots,m$, 
the eigenvalues of the matrix ${A}_k$,
$k=1,\dots,n,\infty$ where the matrix ${A}_\infty$ is defined as
$$
{A}_\infty:= -\sum_{k=1}^{n}{A}_k.
$$
For the sake of technical simplicity let us assume that
\begin{equation}
\lambda^{(k)}_i\neq\lambda^{(k)}_j \,\hbox{for}\quad i\neq j,
\qquad k=1, \dots, n, \infty.
\label{N1.3}
\end{equation}
Moreover, it will be assumed that ${A}_\infty$ is a constant diagonal 
$m\times m$ matrix
with eigenvalues $\lambda^{(\infty)}_j$, $j=1,\dots,m$.

Denote $\Lambda^{(k)}$, $R^{(k)}$ the local monodromy data of the Fuchsian
system near the points $z=u_k$, $k=1$, \dots, $n$, $\infty$. 
The matrices $\Lambda^{(k)}$ are all diagonal
\begin{equation}\label{l-r}
\Lambda^{(k)} = {\rm diag}\, (\lambda_1^{(k)}, \dots, \lambda_m^{(k)}), \quad
k=1, \dots, n, \infty.
\end{equation}
and, under our assumptions
$$
\Lambda^{(\infty)} = {A}_\infty.
$$
Recall that the matrix $G\in GL(m,{\mathbb C})$ belongs to the group ${\mathcal C}_0(\Lambda^{(\infty)})$ {\it iff}
\begin{equation}\label{cinf0}
z^{-\Lambda^{(\infty)}}  G \,  z^{\Lambda^{(\infty)}}
=G_0 +\frac{G_1}{z} +\frac{G_2}{z^2}+\dots .
\end{equation}
It is easy to see that our assumptions about the eigenvalues of $A_\infty$ imply diagonality of the matrix $G_0$.

Let us also remind that the matrices $\Lambda^{(k)}$ satisfy
\begin{equation}\label{l-sum}
\tr\, \Lambda^{(1)} + \dots + \tr\, \Lambda^{(\infty)} = 0.
\end{equation}

\begin{df}
The numbers 
$\lambda_1^{(k)},\dots,\lambda_m^{(k)}$ are called the {\it exponents}\/ of 
the system (\ref{N1}) at the singular point $u_k$.
\end{df}
 
Let us fix a fundamental matrix solutions of the form (\ref{lev6}) near all singular points $u_1$, \dots, $u_n$, $\infty$. To this end
we are to fix branch cuts on the complex plane and choose the branches of
logarithms $\log(z-u_1)$, \dots, $\log (z-u_n)$, $\log z^{-1}$. We will do it 
in the following way: perform parallel branch cuts $\pi_k$ 
between $\infty$ and each of the $u_k$, $k=1,\dots,n$ along a given (generic)
direction. After this we can fix Levelt fundamental matrices analytic on
\begin{equation}\label{domain}
z\in {\mathbb C} \setminus \cup_{k=1}^n \pi_k,
\end{equation}
\begin{equation}
\Phi_k({z})=T_k\left(\ID +{\mathcal O}({z}-u_k)\right) 
({z}-u_k)^{\Lambda^{(k)}}({z}-u_k)^{R^{(k)}}, \quad z\to u_k, \quad 
k=1, \dots, n
\label{N6.1}
\end{equation}
and 
\begin{equation}
\Phi({z})\equiv
\Phi_\infty({z})=\left(\ID+{\mathcal O}({1\over z})\right)
{z}^{-{A}_\infty}{z}^{-{R}^{(\infty)}},\quad\hbox{as}\quad 
{z}\rightarrow\infty,
\label{M3}
\end{equation}
Define the {\it connection matrices}  by 
\begin{equation}
\Phi_\infty({z})=\Phi_k({z}){ C}_k,\label{N3}
\end{equation}
where $\Phi_\infty (z)$ is to be analytically continued in a vicinity of the
pole $u_k$ along the positive side of the branch cut $\pi_k$.

The monodromy matrices ${M}_k$, $k=1,\dots,n,\infty$ are defined
with respect to a basis $l_1,\dots,l_{n}$ of loops in the
fundamental group
$$
\pi_1\left({\mathbb C}\backslash\{u_1,\dots u_{n}\},
\infty\right).
$$ 
Choose the basis in the following way. The loop $l_k$ arrives 
from infinity in a vicinity of 
$u_k$ along one side of the branch cut $\pi_k$ that will be called 
{\it positive}, then it encircles $u_k$ going in anti-clock-wise direction 
leaving all other poles outside and, finally it returns to infinity along 
the opposite side of the branch cut $\pi_k$ called {\it negative}.

Denote $l_j^* \Phi_\infty(z)$ the result of  analytic
continuation of the fundamental matrix $\Phi_\infty(z)$ along the loop
$l_j$. The monodromy matrix ${M}_j$ is defined by
\begin{equation}\label{momo}
l_j^* \Phi_\infty (z) = \Phi_\infty (z) {M}_j, ~~j=1, \dots, n.
\end{equation}
The monodromy matrices satisfy
\begin{equation}
{M}_\infty {M}_{n} \cdots {M}_1=\ID,\qquad
{M}_\infty=\exp\left(2\pi i{A}_\infty\right)
\exp\left(2\pi iR^{(\infty)}\right)
\label{N6}
\end{equation}
if the branch cuts $\pi_1$, \dots, $\pi_{n}$
enter the infinite point according to the order of their labels, i.e., the
positive side of $\pi_{k+1}$  looks at the negative side of $\pi_k$, $k=1,
\dots, n-1$.

Clearly one has
\begin{equation}
{M}_k={ C}_k^{-1} \exp\left(2\pi i \Lambda^{(k)}\right)
\exp\left(2\pi i R^{(k)}\right) 
{ C}_k,\qquad k=1,\dots,n.
\label{N4}
\end{equation}

The collection of the local monodromy data $\Lambda^{(k)}$, $R^{(k)}$ together
with the central connection matrices $C_k$
will be used in order to uniquely fix the Fuchsian system with given poles.
They will be defined up to an equivalence that we now describe. The eigenvalues 
of the diagonal matrices
$\Lambda^{(k)}$
are defined up to permutations. Fixing the order of the eigenvalues, we define
the class of equivalence of the nilpotent part $R^{(k)}$ and of the connection matrices $C_k$ by factoring out the transformations of the form
\begin{eqnarray}\label{class-mon}
&&
R_k \mapsto G_k^{-1} R_k G_k,
\quad
C_k \mapsto G_k ^{-1}C_k G_\infty, \quad k=1, \dots, n,
\nn\\
&&
\nn\\
&& 
G_k \in {\mathcal C}_0 (\Lambda^{(k)}),
\quad G_\infty \in {\mathcal C}_0(\Lambda^{(\infty)}).
\end{eqnarray}
Observe that the monodromy matrices (\ref{N4}) will transform by a simultaneous conjugation
$$
M_k \mapsto G_\infty^{-1} M_k G_\infty, \quad k=1, 2, \dots, n, \infty.
$$
\begin{df} \label{daty}The class of equivalence (\ref{class-mon}) of the collection 
\begin{equation}\label{m-data}
\Lambda^{(1)}, R^{(1)}, \dots, \Lambda^{(\infty)}, R^{(\infty)},  C_1, \dots, C_n
\end{equation}
is called {\it monodromy data} of the Fuchsian system with respect to a fixed ordering of the eigenvalues of the matrices $A_1$, \dots, $A_n$ and a given
choice of the branch cuts.
\end{df}

\begin{lm}
Two Fuchsian systems of the form (\ref{N1}) with the same poles
$u_1,\dots,u_{n},\infty$ and the same matrix $A_\infty$ coincide, modulo diagonal conjugations if and only if they have 
the same monodromy data 
with respect to the same system of branch cuts $\pi_1,\dots,\pi_{n}$. 
\label{lm2.8}
\end{lm}

\begin{proof}
Let 
$$
\Phi_\infty^{(1)}({z})=\left(\ID +O({1\over z})\right)z^{-\Lambda^{(\infty)}}
z^{-R^{(\infty)}}, \quad \Phi_\infty^{(2)}({z})=\left(\ID +O({1\over z})\right)z^{-\tilde\Lambda^{(\infty)}}
z^{-\tilde R^{(\infty)}}
$$ 
be
the fundamental matrices of the form (\ref{M3}) of the two Fuchsian 
systems. Using assumption about $A_\infty$ we derive that $\tilde \Lambda^{(\infty)}=\Lambda^{(\infty)}$.
Multiplying $\Phi_\infty^{(2)}({z})$
if necessary on the right by a matrix 
$G\in {\mathcal C}_0(\Lambda^{(\infty)})$, we can obtain another fundamental
matrix of the second system with 
$$
\tilde R^{(\infty)} = R^{(\infty)}.
$$
Consider the following matrix:
\begin{equation}\label{tYPhi}
Y({z}):= \Phi_\infty^{(2)}({z})[\Phi_\infty^{(1)}({z})]^{-1}.
\end{equation}
$Y({z})$ is an analytic function around infinity: 
\begin{equation}\label{YPhi1}
Y({z})=G_0+{\mathcal O}\left({1\over{z}}\right),\quad\hbox{as}\, 
{z}\rightarrow\infty
\end{equation}
where $G_0$ is a diagonal matrix.
Since the monodromy matrices coincide, $Y({z})$ is a single valued 
function on the punctured Riemann sphere 
$\overline{\mathbb C}\backslash\{u_1,\dots,u_{n}\}$. 
Let us prove that $Y({z})$ is analytic also at the points $u_k$. Indeed, 
having fixed the monodromy data, we can choose the fundamental matrices 
$\Phi_k^{(1)}({z})$ and $\Phi_k^{(2)}({z})$ of the form
(\ref{N6.1}) with the same connection matrices ${ C}_k$
and the same matrices $\Lambda^{(k)}$, $R^{(k)}$. Then near the point
$u_k$, $Y({z})$ is analytic:
\begin{equation}\label{YPhi2}
Y({z})=T_k^{(2)}\left(\ID+{\mathcal O}({z}-u_k) \right)
\left[T_k^{(1)}\left(\ID+{\mathcal O}({z}-u_k) \right)\right]^{-1}.
\end{equation}
This proves that $Y({z})$ is an analytic function on all $\overline{\mathbb C}$ 
and then, by the Liouville theorem
$Y({z})=G_0$, which is constant. So the two Fuchsian systems coincide, after conjugation by the diagonal matrix $G_0$. 
\end{proof}

\begin{rmk} \label{rem1} The connection matrices are determined, within their 
equivalence
classes by the monodromy matrices if the quotients ${\mathcal C}(\Lambda^{(k)},
R^{(k)})/ {\mathcal C}_0(\Lambda^{(k)},
R^{(k)})$ are trivial for all $k=1$, \dots, $n$. In particular this is the case
when all the characteristic exponents at the poles $u_1$, \dots, $u_n$ are
non-resonant.
\end{rmk}

{}From the above Lemma the following result readily follows.

\begin{thm} \label{iso1} If the matrices ${A}_k(u_1, \dots, u_{n})$ satisfy
Schlesinger equations (\ref{sch}) and the matrix 
$$
{A}_\infty = -({A}_1 + \dots + {A}_n)
$$
is diagonal then all the characteristic exponents
do not depend on $u_1$, \dots, $u_{n}$. The fundamental matrix 
$\Phi_\infty(z;u)$ can be chosen in such a way that the nilpotent matrix 
${R}^{(\infty)}$ and also all the monodromy matrices are constant in 
$u_1$, \dots, $u_{n}$. The coefficients of expansion of the fundamental matrix
in $1/z$ belong to a Picard - Vessiot type extension of the field ${\mathcal K}_{(n,m)}$ associated with the solution to Schlesinger equations.
Moreover, the Levelt fundamental matrices $\Phi_k(z;u)$ can be chosen in 
such a way that
all the nilpotent matrices $R^{(k)}$ and also all the connection matrices
${\mathcal C}_k$ are constant. Viceversa, if the deformation 
${A}_k={A}_k(u_1, \dots, u_{n})$ is such that the monodromy data do
not depend on $u_1$, \dots, $u_n$ then the matrices 
${A}_k(u_1, \dots, u_{n})$, $k=1$, \dots, $n$  satisfy Schlesinger
equations.
\end{thm}

Recall that the $u$-dependence of the needed fundamental matrix $\Phi_\infty(z;u)$ is to be determined from the linear equations
\begin{equation}\label{paf}
\partial_i \Phi_\infty(z;u)=-\frac{A_i}{z-u_i}\, \Phi_\infty(z;u), \quad i=1, \dots, n,
\end{equation}
so
\begin{eqnarray}\label{dpsi1}
&&
\partial_i\Psi_1 =-A_i
\\
&&\label{dpsi2}
\partial_i \Psi_2 = -A_i\Psi_1 -u_i A_i
\end{eqnarray}
etc.

\begin{exam}\label{contro}
The following example shows that in general the coefficients of expansion of the fundamental matrix may not be in the field ${\mathcal K}_{(n,m)}$. Indeed,
let us consider the following isomonodromic deformation of the Fuchsian system
$$
\frac{d\Phi}{dz} = \left[ \frac{A_1}{z}
 +\frac{A_2}{z-x}+\frac{A_3}{z-1}\right]\,\Phi,
$$
\begin{eqnarray}
&&
A_1=\left(
\begin{array}{ll}
 -\frac{\left(\sqrt{x}+1\right)^2}{16 \sqrt{x}} & -\frac{1}{2 \sqrt{x}} \\
  & \\
   & \\
 \frac{\left(\sqrt{x}+1\right)^4}{128 \sqrt{x}} & \
\frac{\left(\sqrt{x}+1\right)^2}{16 \sqrt{x}}
\end{array}
\right)
\nn\\
&&
\nn\\
&&
A_2=\left(
\begin{array}{ll}
 -\frac{3 \sqrt{x}-1}{16 \sqrt{x}} & \frac{1}{2 \left(\sqrt{x}+1\right) \
\sqrt{x}} \\
 & \\
 -\frac{\left(\sqrt{x}+1\right) \left(3 \sqrt{x}-1\right)^2}{128 \sqrt{x}} & \
\frac{3 \sqrt{x}-1}{16 \sqrt{x}}
\end{array}
\right)
\nn\\
&&
\nn\\
&&
A_3=\left(
\begin{array}{ll}
 \frac{1}{16} \left(\sqrt{x}-3\right) & \frac{1}{2 \left(\sqrt{x}+1\right)} \
\\
 & \\
 -\frac{1}{128} \left(\sqrt{x}-3\right)^2 \left(\sqrt{x}+1\right) & \
\frac{1}{16} \left(3-\sqrt{x}\right)
\end{array}
\right).
\nn
\end{eqnarray}
In this case
$$
A_\infty = \left(\begin{array}{cc} \frac12 & 0 \\ 0 & -\frac12\end{array}\right), \quad R^{(\infty)}=R^{(\infty)}_1 =\left( \begin{array}{cc}  0 & -\frac12 \\ 0 & 0\end{array}\right).
$$
The fundamental matrix
$$
\Phi =\left( \ID +\frac{\Psi_1}{z}+O\left(\frac1{z^2}\right) \right) \, z^{-A_\infty} z^{-R}
$$
satisfying also the equation
$$
\frac{\partial \Phi}{\partial x} = -\frac{A_2}{z-x} \, \Phi
$$
has
$$
\Psi_1=\left(
\begin{array}{ll}
 \frac{1}{16} \left(3 x-2 \sqrt{x}\right) & -\log \left(\sqrt{x}+1\right) \\
  & \\
 \frac{1}{128} \left(\frac{9 x^2}{2}+2 x^{3/2}-5 x+2 \sqrt{x}\right) & \
\frac{1}{16} \left(2 \sqrt{x}-3 x\right)
\end{array}
\right).
$$
This matrix does not belong to the field ${\mathcal K}_{3,2}$ isomorphic in this case to the field of rational functions in $\sqrt{x}$.
\end{exam}


\setcounter{equation}{0}
\setcounter{theorem}{0}

\section{Reductions of the Schlesinger systems.}

\subsection{Reducible monodromy groups.}

\begin{df}
Given a Fuchsian system of the form (\ref{N1}), we say that its monodromy 
group $\langle {M}_1,\dots,{M}_{n}\rangle$ is 
{\it $l$-reducible,}\/ $0<l< m$ if the monodromy matrices admit a
common invariant subspace $X_l$ of dimension $l$ in the space
of solutions of the system (\ref{N1}).
\end{df}

In particular, if the monodromy group is $l$-reducible, then there
exists a basis where all monodromy matrices have the form
$$
{M}_{k}=\left(\begin{array}{c|c}
\delta_k & \beta_k\\
\hline
0 & \gamma_k
\end{array}\right),\qquad k=1,\dots,n,\infty,
$$
where $\delta_k$, $\beta_k$ and $\gamma_k$ are respectively some $l\times l$, 
${l\times(m-l)}$ and ${(m-l)\times(m-l)}$ matrices.

Given the above definition, we can proceed to the proof of Theorem  \ref{itro1}.

We begin with the proof of Lemma  \ref{perv}.
Our proof, valid for the case of diagonalizable $A_\infty$, is based on the fact that the sum of the exponents of the invariant sub-space $X_l$ must always be a negative integer (see \cite{AB} Lemma 5.2.2). We will perform a sequence of gauge transformations which map such sum to zero. Let $\lambda_1^{(\infty)},\dots,\lambda_m^{(\infty)}$ be the eigenvalues of $A_\infty$ (which is assumed to be diagonal). By means of a permutation $P\in S_m$, we order the eigenvalues of $A_\infty$ 
as follows: the first $l$ eigenvalues correspond to the 
invariant sub-space $X_l$ and we order them in such a way that $\Re{\lambda_1}\geq\Re{\lambda_r}$,
for all $r=2,\dots,l$. Then we order the other eigenvalues in such a way that  $\Re{\lambda_m}\leq\Re{\lambda_s}$ for all $s=l+1,\dots,m-1$. 

Let us fix a fundamental matrix $\Phi$ normalized at infinity
$$
\Phi_\infty= \left(\ID +\frac{\Psi_1}{z} +\frac{\Psi_2}{z^2}+
{\mathcal O}\left(\frac{1}{z^3}\right)\right) z^{- A_\infty} z^{-R^{(\infty)}},
$$
where $\Psi_1$, $\Psi_2$ and $R^{(\infty)}$ are given by formulae (\ref{primer}) in Example \ref{example-psi}.

Consider the following gauge transformation $\Phi(z) = (I(z)+G)\tilde\Phi(z)$  where 
$$
I(z):={\rm Diagonal}\left(z,0,\dots,0\right),
$$
and
\begin{eqnarray}\label{newgauge}
&&
G_{m1}= \Psi_{1_{m1}},\qquad G_{{1m}}=-\frac{1}{G_{m1}},\nn\\
&&
\hbox{if}\quad p\neq 1,m,\quad
G_{{pp}}=1,\quad
G_{{1p}}=\Psi_{1_{mp}} G_{1m},\qquad
G_{{p1}}=\Psi_{1_{p1}},\\
&&
\hbox{if}\quad p,q \neq 1,\, p\neq q,\quad G_{{pq}}=0,\nn\\
&&
G_{{11}}= G_{1m} \Psi_{2_{m1}}+\Psi_{1_{11}}, \quad\hbox{and}\qquad G_{mm}=0.
\nn\end{eqnarray}
Let us first observe that the entries of the matrix $G$ belong to an extension of the differential field ${\mathcal K}_{(n,m)}$ obtained by adding solutions of the linear equations (\ref{dpsi1}), (\ref{dpsi2}).
In order to see that this gauge transformation always works let us show that $\Psi_{1_{m1}}(u)$ is never identically equal to zero if at least one of the $(m, 1)$ matrix entries of the matrices $A_1(u)$, \dots, $A_n(u)$ is different from identical zero. Indeed, this follows from the equations (\ref{dpsi1}).

Let us prove that this transformation maps the matrices $A_1,\dots,A_n$ to new matrices 
$\tilde A_1,\dots,\tilde A_n$ given by
$$
\tilde A_k:= (I(u_k)+G)^{-1} A_k (I(u_k)+G),
$$
such that
\begin{equation}\label{con-gauge}
\tilde A_\infty=-\sum_{k=1}^n \tilde A_k={\rm diagonal}\left(
\lambda^{(\infty)}_1+1,\lambda^{(\infty)}_2,\dots,
\lambda^{(\infty)}_{m-1},\lambda^{(\infty)}_m-1\right).
\end{equation}
In fact $(I(z)+G)^{-1}=J(z)+G^{-1}$ where
$$
J(z):={\rm Diagonal}\left(0,\dots,0,z\right),
$$
therefore
$$
\tilde A_k:=G^{-1} A_k I(u_k) + G^{-1} A_k G + J(u_k) A_k I(u_k)+J(u_k) A_k G.
$$
Multiplying by $G$ from the left and summing on all $k$ we get that the condition (\ref{con-gauge}) is satisfied if and only if
\begin{eqnarray}\label{gauge1}\nn
&&
\left(\begin{array}{ccc}
-g_{11}&(\lambda^{(\infty)}_1-\lambda^{(\infty)}_2) g_{12}&\dots\\
(\lambda^{(\infty)}_2-\lambda^{(\infty)}_1-1) g_{21}&0&\dots\\
\dots&0&\dots\\
(\lambda^{(\infty)}_m-\lambda^{(\infty)}_1-1) g_{m1}&0&\dots\\
\end{array}
 \right.\qquad\qquad\nn\\
 &&
\qquad\qquad\qquad\left.\begin{array}{ccc}
\dots&(\lambda^{(\infty)}_1-\lambda^{(\infty)}_{m-1}) g_{1\,m-1}&
(\lambda^{(\infty)}_1-\lambda^{(\infty)}_m+1) g_{1m}\\
\dots&\dots&0\\
\dots &\dots&0\\
\dots &\dots&0\\
 \end{array}
 \right)=\nn\\
 &&\nn\\
&&\nn
=\left(\begin{array}{cccc}
\sum_{k} A_{k_{11}}u_k&0&\dots&0\\
\dots &0&\dots&0\\
\sum_{k} A_{k_{m1}}u_k&0&\dots&0\\
\end{array}
 \right)+
 \left(\begin{array}{cccc}
g_{1m} \sum_{k} A_{k_{m1}}u_k^2&0&\dots&0\\
0 &0&\dots&0\\
\dots&0&\dots&0\\
\end{array}
 \right)+\\
 &&
 \\
 &&\nn
\qquad\quad +\left(\begin{array}{cccc}
g_{1m}\sum_{s}\sum_{k} A_{k_{ms}}u_k  g_{s1}&\dots&
g_{1m}\sum_{s}\sum_{k} A_{k_{ms}}u_k  g_{sm}\\
0&\dots&0\\
\dots &\dots&\dots\\0&\dots&0
\end{array}
 \right).
 \end{eqnarray}
Observe that in the non-resonant case, these formulae are clearly satisfied thanks to the fact that $\Psi_1$, $\Psi_2$ and $R^{(\infty)}$ are given by formulae (\ref{primer}) in Example \ref{example-psi}. In the resonant case, we only need to prove that when there is a resonance of type  $\lambda^{(\infty)}_m-\lambda^{(\infty)}_p= 1$ or 
$\lambda^{(\infty)}_p-\lambda^{(\infty)}_1= 1$  for any $p=1,\dots,m-1$, then the corresponding coefficients $\sum_k A_{k_{mp}} u_k$, and $\sum_k A_{k_{p1}} u_k$ are zero. Observe that such entries coincide with the $(m,p)$ and $(p,1)$ entries in the matrix $R^{(\infty)}_1$ defined in Section 2.1 (see the formulae (\ref{primer})).
Due to our ordering of the eigenvalues, if $\lambda^{(\infty)}_m-\lambda^{(\infty)}_p= 1$ then $p=1,\dots,l$ and if $\lambda^{(\infty)}_p-\lambda^{(\infty)}_1= 1$ then $p=l+1,\dots,m$. This means that the corresponding $R^{(\infty)}_1$ must lie in the $l\times(m-l)$ lower left block, which is $0$ by the hypothesis that the monodromy group is $l$-reducible. 

Finally, if $\lambda^{(\infty)}_m-\lambda^{(\infty)}_1=2$, we find that 
the gauge transformation works only if 
\begin{eqnarray}\nonumber
&&
\left(\sum_{l=1}^{n}{A}_{l_{m1}}u_l\right)\left(
\sum_{l=1}^{n}\left({A}_{l_{11}}-{A}_{l_{mm}}\right)u_l\right)
-\sum_{l=1}^{n}{A}_{l_{m1}}u_l^2-\nonumber\\
&&
\qquad\qquad -\sum_{p=2}^{m-1}\left(\sum_{l=1}^{n}{A}_{l_{mp}}u_l\right)
G_{{p1}}=0.\nonumber
\end{eqnarray}
This is precisely the condition $\left(R^{(\infty)}_2\right)_{m1}=0$, as it follows from (\ref{primer}). 

Let us prove that this gauge transformation preserves the Schlesinger equations. Differentiating $\tilde A_k$ w.r.t. $u_j$, with $j\neq k$ and using the Schlesinger equations for $A_1,\dots,A_n$ we get:
\begin{eqnarray}\nn
&&
\frac{\partial \tilde A_k}{\partial u_j}=\left[\tilde A_k, (I(u_k)+G)^{-1} \frac{\partial G}{\partial u_j}
+\frac{ (I(u_k)+G)^{-1} A_j(I(u_k)+G)}{u_k-u_j}
\right]=\nn\\
&&
\qquad\qquad=\frac{\left[\tilde A_k, \tilde A_j\right] }{u_k-u_j}+\nn\\
&&
+
 \left[\tilde A_k,(I(u_k)+G)^{-1}\left(\frac{\partial G}{\partial u_j}
+\frac{A_j(I(u_k)-I(u_j))-B_{kj} A_j(I(u_j)+G)
}{u_k-u_j}\right)\right],
\nn\end{eqnarray}
where 
$$
B_{kj}=\left(\begin{array}{cccc}
0&\dots&0&\frac{u_k-u_j}{g_{m1}}\\
0&\dots&0&0\\
\dots&\dots&\dots&\dots\\
0&\dots&0&0
\end{array}\right).
$$
Given the formulae (\ref{newgauge}), it is straightforward to prove that the equation
$$
\frac{\partial G}{\partial u_j}
+\frac{A_j(I(u_k)-I(u_j))-B_{kj} A_j(I(u_j)+G)}{u_k-u_j}=0,
$$
is equivalent to the equations (\ref{dpsi1}), (\ref{dpsi2}).
This proves that also $\tilde A_1,\dots,\tilde A_n$ satisfy the Schlesinger equations.

Now let the sum of the exponents of the invariant sub-space $X_l$ be $-N$, where $N$ is a positive  integer. By iterating the above gauge transformation $N$ times,  we arrive at a new solution  $({B}_1,\dots,{B}_{n})$
of the Schlesinger equations $S_{(n,m)}$ such that the sum  of the exponents of the invariant sub-space $X_l$ is zero and
$$
B_\infty={\rm Diagonal}\left(\lambda^{(\infty)}_1+N,\lambda^{(\infty)}_2,\dots,
\lambda^{(\infty)}_{m-1},\lambda^{(\infty)}_{m}-N\right).
$$
To conclude the proof of this lemma, let us prove that this new  solution  $({B}_1,\dots,{B}_{n})$ is of the form
$$
B_{k_{ij}}=0,\qquad\forall\quad i= l+1,\dots,n,
\quad j=1,\dots,l.
$$
In fact suppose by contradiction that $B_k$ are not in the above form. Then by Lemma 5.2.2. in \cite{AB}, there exists a gauge transformation $P$, constant in $z$, such that
the new residue matrices $\tilde{B}_k=P^{-1}{B}_kP$ have the form
has the form 
$$
\tilde{B}_{k_{ij}}=0,\qquad\forall\quad i= l+1,\dots,n,
\quad j=1,\dots,l.
$$
In general $\tilde B_\infty$ won't be diagonal, but we can diagonalize it by a constant gauge transformation $Q$ preserving the block triangular form of $\tilde B_1,\dots,\tilde B_n$. So we end up with
$$
\hat B_\infty=Q^{-1} P^{-1} B_\infty P Q,\qquad \hat B_k = Q^{-1} P^{-1} B_k P Q,
$$
and since $Q^{-1} P^{-1} B_\infty P Q=B_\infty$, we have that $P Q$ is diagonal. But then if $B_k$ is not block triangular, $\hat B_k$ is not either, so we obtain a contradiction. Lemma \ref{perv} is proved.
{$\quad$\hfill
\raisebox{0.11truecm}{\fbox{}}\par\vskip0.4truecm}

\vskip 0.2 cm
\noindent{\bf Proof of Theorem \ref{itro1}}. 
By Lemma \ref{perv}, we obtained a gauge transformation mapping 
a solution $({A}_1,\dots,{A}_{n})$ 
of the Schlesinger system ${S}_{n,m}$ with a $l$-reducible monodromy group 
to a solution $B_1,\dots,B_n$ of the block triangular form.
As it was explained in the Introduction, the solution
$({B}_1(u),\dots,{B}_{n}(u))$ belongs to a Picard--Vessiot type extension 
${\mathcal K}^{(N)}$ for some $N$
of the composite 
$$
{\mathcal K} = {\mathcal K}_{n,l} {\mathcal K}_{n,m-l}.
$$
So, to conclude the proof of this theorem, we need to prove Lemma  \ref{vtor}.

Let us prove the formulae (\ref{sviaz}). Our gauge transformation constructed in Lemma \ref{perv}
 is an iteration of elementary gauges transformation $\Phi=(I(z)+G)\tilde\Phi$
mapping the matrices $A_1,\dots,A_n$ to new matrices $\tilde A_1,\dots,\tilde A_n$ such that
$\tilde A_\infty=A_\infty+{\rm Diagonal}(1,0,\dots,0,-1)$. 

Let us prove that each elementary gauge transformation preserves the normalization at infinity.
More precisely, we prove that if we fix a fundamental matrix $\Phi$ normalized at infinity
$$
\Phi_\infty= \left(\ID +\frac{\Psi_1}{z} +\frac{\Psi_2}{z^2}+
{\mathcal O}\left(\frac{1}{z^3}\right)\right) z^{- A_\infty} z^{-R^{(\infty)}},
$$
then $\tilde\Phi=(J(z)+G^{-1})\Phi_\infty= \left(\ID +
{\mathcal O}\left(\frac{1}{z}\right)\right) z^{-\tilde A_\infty} z^{-\tilde R^{(\infty)}}$, whith $\tilde R^{(\infty)}= R^{(\infty)}$.

In fact it is straightforward to prove that
\begin{eqnarray}\nn
&&
(J(z)+G^{-1}) \left(\ID +\frac{\Psi_1}{z} +\frac{\Psi_2}{z^2}+
{\mathcal O}\left(\frac{1}{z^3}\right)\right){\rm Diagonal}(z,0,\dots,0,\frac{1}{z}) =\nn\\
&&
\qquad\qquad=\chi_1 z + \chi_0 +{\mathcal O}\left(\frac{1}{z}\right)\nn
\end{eqnarray}
where all matrix elements of $\chi_1$ are zero apart from the $(m,1)$ element which is
$$
\chi_{1_{m1}}= \frac{1}{g_{1m}}+\Psi_{1_{m1}}
$$
and the matrix elements of $\chi_0$ are given by the following: for $p\neq 1,m$
$$
\chi_{0_{pp}}=1, \quad
\chi_{0_{p1}} = -g_{p1}+\Psi_{1_{p1}},\quad
\chi_{0_{pm}}=  -\frac{g_{1p}}{g_{1m}}+  \Psi_{1_{mp}} ,
$$
$$
\chi_{0_{11}}=\chi_{0_{mm}}=1,
$$
and 
$$\chi_{0_{m1}}  =  -\frac{g_{11} -\sum_{p=2}^{m-1} g_{1p} \Psi_{1_{pm}} -\Psi_{1_{11}}+ \sum_{p=2}^{m-1}  g_{1p} \Psi_{1_{p1}}}{g_{1m}}+\Psi_{2_{1m}}.
$$
Using the formulae (\ref{newgauge}) for $G$ it is easy to prove that all entries of $\chi_0$ and $\chi_1$ are zero.

Therefore each elementary gauge transformation preserves the normalization at infinity and maps $A_\infty$ to 
$$
\tilde A_\infty=A_\infty+{\rm Diagonal}(1,0,\dots,0,-1).
$$

Since the fundamental matrix remains normalized at infinity and the gauge transformation $\Phi=(I(z)+G)\tilde\Phi$ is analytic over $\mathbb C$, all monodromy data
$\Lambda^{(1)}(A), R^{(1)}(A), \dots, \Lambda^{(n)}(A), R^{(n)}(A),  C_1(A), \dots, C_n(A)$
are preserved in each iteration. Finally we prove that $R^{(\infty)} =\tilde R^{(\infty)}$. Due to the above we only need to prove that if
$$
z^{-A_\infty} R^{(\infty)}(A) z^{A_\infty}= \frac{R_1}{z} + \frac{R_2}{z^2} + ......
$$
where $R_1,R_2,\dots$ are some matrices defined in Section 2, then the matrix
$$
z^{-B_\infty} R^{(\infty)}(A) z^{B_\infty}
$$
is also polynomial in $\frac{1}{z}$.
Since $B_\infty=A_\infty +{\rm Diagonal}(N,0,\dots,0,-N)$ we get
\begin{eqnarray}\nn
&&
z^{-B_\infty} R^{(\infty)}(A) z^{B_\infty} =\nn\\
&&
{\rm Diagonal}(z^{-N},1,\dots,1,z^N)
z^{-A_\infty} R^{(\infty)}(A) z^{A_\infty}
{\rm Diagonal}(z^{N},1,\dots,1,z^{-N}) =\nn\\
&&
{\rm Diagonal}(z^{-N},1,\dots,1,z^N)
\left(\frac{R_1}{z} + \frac{R_2}{z^2} + ......\right)
{\rm Diagonal}(z^{N},1,\dots,1,z^{-N})=\nn\\
&&
Pol\left(\frac{1}{z}\right) + Div(z),
\nn\end{eqnarray}
where $Pol\left(\frac{1}{z}\right) $ and $Div(z)$ are matrix values  polynomials in $\frac{1}{z}$  and $z$ respectively.  The matrix elements of the latter are of the form:
$$
Div_{pq}=0,\quad\hbox{for}\quad q\neq 1, p\neq m,\quad
Div_{11}=0,\quad Div_{mm}=0,
$$
$$
Div_{p1} = \sum R^{(\infty)}_{k_{p1}} z^{N-k},\quad\hbox{for}\quad p\neq 1,m,
$$
$$
Div_{mq} = \sum R^{(\infty)}_{k_{mq}} z^{N-k},\quad\hbox{for}\quad q\neq 1,m
$$
$$
Div_{m1} = \sum R^{(\infty)}_{k_{m1}} z^{2N-k}.
$$
Since the monodromy group is reducible, all the entries of $R^{(\infty)}_k$ involved in the above expressions are identically zero. Therefore $Div(z)\equiv 0$ as we wanted to prove. 
This proves the relations (\ref{sviaz}).

Let us now prove the statement ii) of Lemma \ref{vtor}.
Starting form the 
solution $(B_1,\dots,B_{n})$, we can reconstruct $({A}_1,\dots,{A}_{n})$ by iterating another gauge transformation of the form $(J(z)+F)$ where $J(z)={\rm Diagonal}(0,\dots,0,z)$ and
\begin{eqnarray}\label{new-inv}
&&
F_{1m}=\tilde\Psi_{1_{1m}},\qquad 
F_{m1}=-\frac{1}{\tilde\Psi_{1_{1m}}},\nn\\
&& 
\hbox{if}\quad p\neq 1,m,\quad
F_{{pp}}=1,\quad
F_{{mp}}=\tilde\Psi_{1_{1p}} G_{m1},\qquad
F_{{pm}}=\tilde\Psi_{1_{pm}},\\
&&
\hbox{if}\quad p,q \neq 1,\, p\neq q,\quad F_{{pq}}=0,\nn\\
&&
F_{{11}}=0, \quad\hbox{and}\qquad 
F_{mm}=F_{m1} \Psi_{2_{1m}}+\Psi_{1_{mm}}.
\nn\end{eqnarray}
This gauge transformation is always well defined because $\tilde\Psi_{1_{1m}}$ is always non-zero (proof of this fact is analogous to the proof that $\Psi_{1_{m1}}$ is never zero given above).
Following the same computations as in the proof of  Lemmata \ref{perv} and \ref{vtor}, it is easy to verify that this gauge transformation preserves the Schlesinger equations, the normalization of the fundamental matrix at infinity, $R^{(\infty)}$ and maps $\tilde A_\infty$ to 
$$
A_\infty=\tilde A_\infty-{\rm Diagonal}(1,0,\dots,0,-1).
$$
The above arguments complete the proof of Lemma \ref{vtor} and, therefore of Theorem \ref{itro1}. 
{$\quad$\hfill
\raisebox{0.11truecm}{\fbox{}}\par\vskip0.4truecm}

\subsubsection{Upper triangular monodromy groups.}

In this section we deal with the case of upper triangular monodromy
groups, that is there exists a basis where all monodromy
matrices have the form
$$
{M}_{k_{ij}}=0,\qquad\forall\qquad i>j.
$$

To prove Corollary \ref{tria}
we iterate the procedure of the proof of Lemma \ref{perv}: at the 
first step we show that $({A}_1,\dots,{A}_{n})$ is mapped by a rational gauge
transformation to $({A}_1^{(1)},\dots,{A}_{n}^{(1)})$ of the form
$$
{{A}_k^{(1)}}_{i1}=0,\qquad \forall\,i\neq 1, \quad \forall\,k=1,\dots,n.
$$
At the $l$-th step we show that $({A}_1,\dots,{A}_{n})$ is mapped by a 
rational gauge transformation to $({A}_1^{(l)},\dots,{A}_{n}^{(l)})$ of 
the form
$$
{{A}_k^{(l)}}_{ij}=0,\qquad i>j, j\leq l, \quad k=1,\dots,n.
$$
At the $m$-th step we obtain that is mapped by a 
rational gauge transformation to 
$$
\tilde{A}_{k_{ij}}:={{A}_{k_{ij}}^{(m)}}=0,\qquad\forall\qquad i>j.
$$
Let us show that $\tilde{A}_{k_{ij}}(u_1,\dots,u_{n})$ belongs to the 
Picard--Vessiot type extension 
${\mathcal K}^{(N)}$ for some $N$ of
$$
{\mathcal K}={\mathbb C}(u_1,\dots,u_n).
$$
Clearly the diagonal
elements $\tilde{A}_{k}$ are the eigenvalues 
$\lambda^{(k)}_1,\dots \lambda^{(k)}_m$.
The Schlesinger equations for $i\neq j$ read:
\begin{eqnarray}\nonumber
&&
{\partial\over\partial u_j} \tilde{A}_{i_{p,p+q}}{\lambda^{(j)}_{p+q}-\lambda^{(j)}_{p}\over u_i-u_j}\tilde{A}_{i_{p,p+q}}-
{\lambda^{(i)}_{p+q}-\lambda^{(i)}_{p}\over u_i-u_j}
\tilde{A}_{j_{p,p+q}}+\nonumber\\
&&
+\sum_{s=1}^{q-1}{\tilde{A}_{i_{p,p+s}}\tilde{A}_{j_{p+s,p+q}}-
\tilde{A}_{j_{p,p+s}}\tilde{A}_{i_{p+s,p+q}}\over u_i-u_j},
\end{eqnarray}
and for $i=j$
\begin{eqnarray}\nonumber
&&
{\partial\over\partial u_i} \tilde{A}_{i_{p,p+q}}= 
-\sum_{j\neq i}\big[
{\lambda^{(j)}_{p+q}-\lambda^{(j)}_{p}\over u_i-u_j}
\tilde{A}_{i_{p,p+q}}-
{\lambda^{(i)}_{p+q}-\lambda^{(i)}_{p}\over u_i-u_j}
\tilde{A}_{j_{p,p+q}}+\nonumber\\
&&
+\sum_{s=1}^{q-1}{\tilde{A}_{i_{p,p+s}}
\tilde{A}_{j_{p+s,p+q}}-
\tilde{A}_{j_{p,p+s}}\tilde{A}_{i_{p+s,p+q}}
\over u_i-u_j}\big],\nonumber
\end{eqnarray}
where for $q=1$ the sum $\sum_{s=1}^{q-1}$ is zero. 
It is clear that for each $q$, $1\leq q< m-p$, the differential
system for the matrix elements ${A}_{i_{p,p+q}}$, $i=1,\dots,n$ is 
linear and it is Pfaffian integrable because the Schlesinger equations are 
Pfaffian integrable. In particular it is worth observing that
for each $q$, $1\leq q< m-p$, the homogeneous part of such differential
system is the Lauricella hypergeometric system (see \cite{IKSY}). 



\subsection{An example.}

Consider the following solution $A_1,A_2,A_3$ of the Schlesinger equations in dimension $m=2$, where we have chosen $u_1=0, u_2=x,u_3=1$, with matrix entries:
$$
A_{1_{11}}=
\frac{2 \log{\frac{\sqrt{x}+1}{\sqrt{x}-1}}\sqrt{x}(x^2+4x-5)-4 x(4+3x)-
\left(\log{\frac{\sqrt{x}+1}{\sqrt{x}-1}}\right)^2(x-1)^2}
{2\left((x-1)  \log{\frac{\sqrt{x}+1}{\sqrt{x}-1}}-2\sqrt{x}\right)^2},
$$
\begin{eqnarray}\nn
&&
A_{1_{12}}=
\frac{\left(\log{\frac{\sqrt{x}+1}{\sqrt{x}-1}}\right)^2(x-1)^2+2 x(5+3x)-(x^2+6 x -7)\sqrt{x}
\log{\frac{\sqrt{x}+1}{\sqrt{x}-1}}}
{4\left((x-1)  \log{\frac{\sqrt{x}+1}{\sqrt{x}-1}}-2\sqrt{x}\right)^2(x-1)}\times\nn\\
&&
\qquad\qquad\times
\left(
6\sqrt{x}(1+x)-(x^2+2 x-3) \log{\frac{\sqrt{x}+1}{\sqrt{x}-1}}
\right)\nn\end{eqnarray}
$$
A_{1_{21}}=\frac{4\sqrt{x}(1-x)}{\left((x-1)  \log{\frac{\sqrt{x}+1}{\sqrt{x}-1}}-2\sqrt{x}\right)^2}\nn\\
$$
$$
A_{2_{11}}=\frac{
\left(\log{\frac{\sqrt{x}+1}{\sqrt{x}-1}}\right)^2(x-1)^3-4 x(7+x)-
8  \log{\frac{\sqrt{x}+1}{\sqrt{x}-1}}\sqrt{x}(x^2-3x+2)}
{4\left((x-1)  \log{\frac{\sqrt{x}+1}{\sqrt{x}-1}}-2\sqrt{x}\right)^2(x-1)},
$$
\begin{eqnarray}\nn
&&
A_{2_{12}}=
-\frac{\left(\log{\frac{\sqrt{x}+1}{\sqrt{x}-1}}\right)^2(x-1)^3-16 x-2 (3x^2-8 x+5)\sqrt{x}
\log{\frac{\sqrt{x}+1}{\sqrt{x}-1}}}
{8\left((x-1)  \log{\frac{\sqrt{x}+1}{\sqrt{x}-1}}-2\sqrt{x}\right)^2(x-1)^2}\times\nn\\
&&
\qquad\qquad\times
\left(
2\sqrt{x}(3+x)+(x^2-4 x+3) \log{\frac{\sqrt{x}+1}{\sqrt{x}-1}}
\right)\nn\end{eqnarray}
$$
A_{2_{21}}=-\frac{4\sqrt{x}}{\left((x-1)  \log{\frac{\sqrt{x}+1}{\sqrt{x}-1}}-2\sqrt{x}\right)^2}
$$
$$
A_{3_{11}}=\frac{1}{2}-A_{1_{11}}-A_{2_{11}},\qquad
A_{3_{12}}=-A_{1_{12}}-A_{2_{12}},\qquad
A_{3_{21}}=-A_{1_{21}}-A_{2_{21}},
$$
$$
A_{1_{22}}=-A_{1_{11}},\qquad A_{2_{22}}=-A_{2_{11}},\qquad A_{3_{22}}=-A_{3_{11}}.
$$
This solution has a reducible monodromy group.
Observe that 
$$
A_\infty=\left(\begin{array}{cc}-\frac{1}{2}&0\\0&\frac{1}{2}\end{array} \right)
$$ 
is resonant. By applying our technique, it is straightforward to obtain a the new solution $B_1,B_2,B_3$ of the Schlesinger equations, gauge equivalent to $A_1,A_2,A_3$ in the upper triangular form:
$$
B_1=\left(\begin{array}{cc}
\frac{1}{2}&\frac{\sqrt{x}}{x-1}\\
&\\
0&-\frac{1}{2}
\end{array}\right),\qquad
B_2=\left(\begin{array}{cc}
-\frac{1}{4}&\frac{\sqrt{x}}{(x-1)^2}\\
&\\
0&\frac{1}{4}
\end{array}\right),
$$
and
$$
B_3=\left(\begin{array}{cc}
-\frac{3}{4}&-\frac{x\sqrt{x}}{(x-1)^2}\\
&\\
0&\frac{3}{4}
\end{array}\right),\qquad B_\infty=\left(\begin{array}{cc}
\frac{1}{2}&0\\ 0&-\frac{1}{2}\\
\end{array}\right).
$$
This new solution is actually algebraic. 
This shows that  the differential fields ${\mathcal K}_{3,2}^A$ and ${\mathcal K}_{3,2}^B$ associated with the solutions $A_1,A_2,A_3$ and $B_1,B_2,B_3$ respectively are not isomorphic.

\subsection{Smaller monodromy groups.}
The proof of Theorem \ref{l-malo} is based on a few lemmata:

\begin{lm}
Let ${A}_{1},\dots,{A}_{n}$ be a solution of the Schlesinger 
equations such that one of the monodromy matrices
$({M}_1,\dots,{M}_{n})$, say  
${M}_l$, is proportional to the identity, 
then then there exists a solution 
$\tilde {A}_1,$\dots,$\tilde {A}_{l-1},\tilde {A}_{l+1},$
\dots$\tilde {A}_{n}$ of the Schlesinger equations in $n-1$
variables with monodromy matrices
${M}_1,\dots,{M}_{l-1}$, ${M}_{l+1},
\dots,{M}_{n}$. The original solution ${A}_{1},\dots,{A}_{n}$ depends 
rationally on $\tilde {A}_1$,\dots, $\tilde {A}_{l-1}$,
$\tilde {A}_{l+1},$\dots$\tilde {A}_{n}$, $\tilde\Phi(u_l)$ and
on $u_l$.
\label{small2}\end{lm}

\begin{proof} 
Let us first consider the case ${M}_l=\ID$, for 
simplicity, $l=n$.
This means that all eigenvalues $\lambda^{(n)}_1,\dots \lambda^{(n)}_m$ 
of ${A}_n$ are integers and $R^{(n)}=0$. To eliminate the singularity $n$, we perform a confromal transformation $\zeta=\frac{1}{z-u_n}$. We obtain
$$
\frac{{\rm d}\Phi}{{\rm d}\zeta} = \left(\frac{A_\infty}{\zeta} +
\sum_{k=1}^{n-1} \frac{A_k}{\zeta-\tilde u_k}\right)\Phi,
$$
where $\tilde u_k=\frac{1}{u_k-u_n}$, for $k\neq n$. The new residue matrix at infinity is $-A_n$. We perform a gauge transformation diagonalizing $A_n$ and 
use iterations of the gauge transformation of the form $(I(\zeta)+G)$ where $G$ is defined by 
formulae (\ref{gauge1}) to map all eigenvalues of $A_n$ to zero. 

We have seen in the proof of Lemma \ref{perv} that this gauge transformation is always well defined and it works for $R^{(\infty)}=0$. Of course similar formulae can be given to map any $\lambda^{(\infty)}_j$ to $\lambda^{(\infty)}_j+1$ and any  
$\lambda^{(\infty)}_i$ to 
$\lambda^{(\infty)}_i-1$. After enough iterations we end up with a new Fuchsian system of the form 
$$
\frac{{\rm d}\tilde\Phi}{{\rm d}\zeta} = \left(\frac{\tilde A_\infty}{\zeta}+
\sum_{k=1}^{n-1} \frac{\tilde A_k}{\zeta-\tilde u_k}\right)\tilde\Phi,
$$
such that the residue at infinity is $\tilde A_n=0$.

Now we perform the inverse conformal transformation, $z=\frac{1}{\zeta}+u_n$, we obtain
$$
\frac{{\rm d}\tilde\Phi}{{\rm d}z} = \sum_{k=1}^{n-1} \frac{\tilde A_k}{z- u_k}\tilde\Phi,
$$
and the residue at infinity is $\tilde A_\infty$. We finally perform a gauge transformation diagonalizing $\tilde A_\infty$, so that the final Fuchsian system is
$$
\frac{{\rm d}\hat\Phi}{{\rm d}z} = \sum_{k=1}^{n-1} \frac{\hat A_k}{z- u_k}\hat\Phi,
$$
where $\hat A_\infty=A_\infty$. 

All the monodromy data of this new system coincide with the ones of the original system with matrices ${A}_{1},\dots,{A}_{n},A_\infty$. The proof of this fact is very similar to the proof of statement ii of lemma
\ref{vtor} and we omit it.

The new matrices $\hat A_1,\dots,\hat A_{n-1}$  satisfy the Schlesinger equations because the gauge transformations of the form $(I(z)+G)$ where $G$ is defined by the formulae (\ref{gauge1}) preserve the Schlesinger equations. Observe that since $\hat A_n$ is zero, $\hat A_1,\dots,\hat A_{n-1}$  satisfy the Schlesinger equations $S_{n-1,m}$.

We now want to reconstruct the original solution 
${A}_{1},\dots,{A}_{n}$ from 
$\hat {A}_1,$\dots,$\hat {A}_{n-1}$. 

Let us consider the Fuchsian system
$$
{{\rm d}\hat\Phi\over{\rm d}z}=\sum_{k=1}^{n-1}{\hat {A}_{k}
\over z-u_k}\hat\Phi.
$$
Let us choose any point $u_n\neq u_k$, $k=1,\dots,n-1$ and perform 
the constant gauge transformation
$\hat\Psi=\hat\Phi(u_n)^{-1}\hat\Phi$, 
where $\hat\Phi(u_n)$ is the value at $z=u_n$ of
$$
\hat\Phi(z)\left( \ID +{\mathcal O}({1\over z})\right)
z^{-{A}_\infty} z^{R^{(\infty)}}.
$$
Let us perform the conformal transformation $\zeta={1\over z-u_n}$,
$$
\hat\Psi(\zeta):=\hat\Phi(u_n)^{-1}\hat\Phi({1\over\zeta}+u_n).
$$
Let us apply a product $F_\infty(\zeta)$ of gauge transformations of the form  $(J(\zeta)+F)$,
where $F$ is given by the formulae (\ref{new-inv}), to create a new non-zero residue matrix at 
infinity with integer entries 
$\lambda^{(n)}_1,\dots \lambda^{(n)}_m$:
$$
\hat\Psi(\zeta)=F_\infty(\zeta)\hat\Psi(\zeta)=F_\infty(\zeta)
\hat\Phi(u_n)^{-1}\hat\Phi({1\over\zeta}+u_n).
$$
Let us now apply the conformal transformation $z={1\over\zeta}+u_n$:
$$
\tilde\Phi(z):= 
F_\infty({1\over z-u_n})\hat\Phi(u_n)^{-1}\hat\Phi(z).
$$
We need now to diagonalize the new residue matrix at infinity 
$$
\tilde A_\infty= F_\infty(u_n)\hat\Phi(u_n)^{-1}\hat A_\infty\hat\Phi(u_n)
F_\infty(u_n)^{-1}. 
$$
To do so we put
$$
\Phi(z):= \hat\Phi(u_n)F_\infty(u_n)^{-1}\hat\Phi(z).
$$
The new residue matrices are
$$
{B}_i= \hat\Phi(u_n)F_\infty(u_n)^{-1} 
F_\infty({1\over u_i-u_n})\hat\Phi(u_n)^{-1}
\hat A_i \hat\Phi(u_n)
F_\infty^{-1}({1\over u_i-u_n})F_\infty(u_n)\hat\Phi(u_n)^{-1}
$$
for $i=1,\dots,n-1$ and 
$${B}_n= \hat\Phi(u_n)F_\infty^{-1}(u_n){\rm diagonal}
(\lambda^{(n)}_1,\dots \lambda^{(n)}_m)
F_\infty(u_n)\hat\Phi(u_n)^{-1}.
$$
The Fuchsian system with residue matrices $B_1,\dots,B_n,B_\infty$ has the same exponents and the same monodromy data as the original system of residue matrices $A_1,\dots,A_n,A_\infty$. Therefore,
by the uniqueness lemma \ref{lm2.8},  $A_1,\dots,A_n,A_\infty$ coincide with
 $B_1,\dots,B_n,B_\infty$ up to diagonal conjugation.

As a consequence ${A}_1,$\dots,${A}_{n}$
depend rationally on $\hat {A}_1,$\dots,$\hat {A}_{n-1}$, on  
$\hat\Phi(u_n)$ and $u_n$.

Now let us suppose that ${M}_l$ is only proportional to the identity.
This means that all eigenvalues $\lambda^{(l)}_1,\dots \lambda^{(l)}_m$ of 
${A}_l$ are resonant. Since their sum is zero, the only possibility is
${M}_l=\exp\left({2\pi i s\over m}\right)\ID$ for some $s=1,\dots,m-1$.
To transform this matrix to the identity we use iterations of the symmetries (\ref{simmk}), (\ref{simminf}), to  map our 
solution ${A}_1,\dots,{A}_{n}$ to a solution 
$\hat {A}_1,\dots,\hat {A}_{n}$ having ${M}_n=\ID$. Since these symmetries are birational,  ${A}_1,\dots,{A}_{n}$  are rational functions of $\hat {A}_1,\dots,\hat {A}_{n}$. 
Then we can apply the above procedure to kill $\hat A_n$.

\begin{rmk}
Observe that in Lemma \ref{small2}, for $M_l=\exp\left({2\pi\over m} \right)\ID$, the new solution
$\tilde {A}_1,$\dots,$\tilde {A}_{l-1},\tilde {A}_{l+1},$\dots$\tilde {A}_{n}$ 
has monodromy matrices ${M}_1,\dots,{M}_{l-1},{M}_{l+1}$,
$\dots,{M}_{n}$, and a new monodromy matrix at infinity $\exp\left(-{2\pi\over m} \right)M_\infty$.
\end{rmk}

\begin{lm}
Let $({A}_{1},\dots{A}_{n})$ be a solution of the 
Schlesinger equations with ${M}_{\infty}$ proportional to the identity $\ID$,
say  ${M}_{\infty}=\exp({2\pi i\over m})\ID$.
Suppose that ${M}_n$ is not proportional to the identity, then there exists 
a solution $\tilde {A}_1$,\dots$\tilde {A}_{n-1}$ of the Schlesinger equations 
with monodromy matrices 
\begin{equation}
{\mathcal C}_n{M}_1{\mathcal C}_n^{-1},
\dots,{\mathcal C}_n{M}_{n-1}{\mathcal C}_n^{-1},
\label{monnew}\end{equation}
and $\tilde M_\infty = {\mathcal C}_n\exp(-{2\pi i\over m})M_n{\mathcal C}_n^{-1}$,
${\mathcal C}_n$ being the connection matrix of ${M}_n$.
The given solution ${A}_{1},\dots{A}_{n}$ depends rationally 
on $\tilde {A}_1$,\dots,$\tilde {A}_{n-1}$, $\tilde\Phi(u_n)$ and on $u_n$.
\label{small4}\end{lm}

We perform a symmetry (\ref{simminf}) (or a conformal transformation), in order to apply Lemma \ref{small2}
to the case $M_1$ proportional to the identity.

\noindent{\bf End of the proof of Theorem \ref{small1}.} Suppose that
${A}_{1},\dots,{A}_{n}$ is a solution of the Schlesinger 
equations such that the collection of its monodromy matrices 
${M}_1$,\dots,${M}_{n}$, ${M}_\infty$ is $l$-smaller. 
If none of the monodromy matrices being proportional to the identity is
equal to ${M}_\infty$, we can simply conclude
by $l$ iterations of Lemma \ref{small2}. If ${M}_\infty$ is proportional to 
the identity, first we apply Lemma \ref{small4}, then we iterate Lemma 
\ref{small2} $l-1$ times. This concludes the proof of
Theorem \ref{small1}. 
\end{proof}

\bibliography{bibgenerale1}

\begin{thebibliography}{10}

\bibitem{kitaev}
F.V. Andreev and A.V. Kitaev.
\newblock Transformations \text{$RS\sp 2\sb 4(3)$} of the ranks $\leq 4$ and
  algebraic solutions of the sixth \text{Painlev\'e} equation.
\newblock {\em Comm. Math. Phys.}, \textbf{228}:no. 1, 151--176, 2002.

\bibitem{AB}
D.V. Anosov and A.A. Bolibruch.
\newblock {\em The \text{Riemann-Hilbert Problem}}, volume \textbf{E 22}.
\newblock Aspects of Mathematics, 1994.

\bibitem{Bol2}
A.A. Bolibruch.
\newblock {\em The $21$-st \text{H}ilbert problem for linear \text{Fuchsian}
  systems}.
\newblock Developments in mathematics: the Moscow school. Chapman and Hall,
  London, 1993.

\bibitem{Bol}
A.A. Bolibruch.
\newblock \text{On isomonodromic deformations of Fuchsian systems}.
\newblock {\em J. Dynam. Control Systems}, \textbf{3}:no. 4, 589--604, 1997.

\bibitem{costin}
O.~Costin and R.~D. Costin.
\newblock Asymptotic properties of a family of solutions of the
  \text{P}ainlev\text{\'e} equation \text{VI}.
\newblock {\em Int. Math. Res. Not.}, \textbf{2002}:no.22:1167--1182, 2002.

\bibitem{Dub1}
B.~Dubrovin.
\newblock \text{Painlev\'e} transcendents in \text{Two-Dimensional}
  \text{Topological} \text{Field Theory}.
\newblock {\em The Painlev\'e property One Century Later}, pages 287--412,
  1999.

\bibitem{DM}
B.~Dubrovin and M.~Mazzocco.
\newblock \text{Monodromy of certain Painlev\'e-VI} transcendents and
  reflection groups.
\newblock {\em Invent. Math.}, \textbf{141}:55--147, 2000.

\bibitem{DM1}
B.~Dubrovin and M.~Mazzocco.
\newblock Canonical structures of the \text{Schlesinger} systems.
\newblock {\em preprint math.DG/0311261}, 2003.

\bibitem{fuchs}
R.~Fuchs.
\newblock \text{Ueber} lineare homogene \text{Differentialgleichungen} zweiter
  \text{Ordnung} mit drei im \text{Endlichen} gelegenen wesentlich
  \text{singul\"aren} \text{Stellen}.
\newblock {\em Math. Ann.}, \textbf{63}:301--321, 1907.

\bibitem{Gar1}
R.~Garnier.
\newblock Sur des \'equations diff\'erentielles du troisi\`eme ordre dont
  l'int\'egrale g\'en\'erale est uniforme et sur une classe d'\'equations
  nouvelles d'ordre sup\'erieur dont l'int\'egrale g\'en\'erale a ses points
  critiques fixes.
\newblock {\em Ann. Sci. \'Ecole Norm. Sup.}, \textbf{29}:no. 3, 1--126, 1912.

\bibitem{Gar2}
R.~Garnier.
\newblock Solution du probleme de \text{Riemann} pour les systemes
  diff\'erentielles lin\'eaires du second ordre.
\newblock {\em Ann. Sci. \'Ecole Norm.. Sup.}, \textbf{43}:239--352, 1926.

\bibitem{guzzetti}
D.~Guzzetti.
\newblock The elliptic representation of the general \text{Painlev\'e VI}
  equation.
\newblock {\em Comm. Pure Appl. Math.}, \textbf{55}:no. 10, 1280--1363, 2002.

\bibitem{Hit}
N.~Hitchin.
\newblock A lecture on the octahedron.
\newblock {\em Bull. London Math. Soc.}, \textbf{35}:no.5, 577--600, 2003.

\bibitem{ince}
E.L. Ince.
\newblock {\em Ordinary differential equations}.
\newblock Dover Publications INC., 1956.

\bibitem{its}
A.R. Its and V.Yu. Novokshenov.
\newblock {\em The isomonodromic deformation method in the theory of
  \text{Painlev\'e} equations}, volume \textbf{1191} of {\em Lecture notes in
  mathematics}.
\newblock Springer, 1980.

\bibitem{IKSY}
K.~Iwasaki, H.~Kimura, S.~Shimomura, and M.~Yoshida.
\newblock {\em From \text{Gauss to Painlev\'e}, a Modern Theory of Special
  Functions}, volume \textbf{ E 16}.
\newblock Aspects of Mathematics, 1991.

\bibitem{jimbo}
M.~Jimbo.
\newblock Monodromy problem and the boundary condition for some
  \text{P}ainlev{\'e} equations.
\newblock {\em Publ. Res. Inst. Math. Sci.}, \textbf{18}:1137--1161, 1982.

\bibitem{MJ1}
M.~Jimbo and T.~Miwa.
\newblock Monodromy preserving deformations of linear ordinary differential
  equations with rational coefficients \text{II}.
\newblock {\em Physica 2D}, \textbf{2}:no. 3, 407--448, 1981.

\bibitem{MJU}
M.~Jimbo, T.~Miwa, and K.~Ueno.
\newblock Monodromy preserving deformations of linear ordinary differential
  equations with rational coefficients \text{I}.
\newblock {\em Physica 2D}, \textbf{2}:no. 2, 306--352, 1981.

\bibitem{kol}
E.~R. Kolchin.
\newblock {\em Differential algebra and algebraic groups}.
\newblock Academic Press, New York, 1973.
\newblock Pure and Applied Mathematics, Vol. 54.

\bibitem{Levelt}
A.H.M. Levelt.
\newblock Hypergeometric functions.
\newblock {\em Doctoral thesis, University of Amsterdam}, 1961.

\bibitem{malek1}
S.~Malek.
\newblock \text{A Fuchsian} system with a reducible monodromy is
  meromorphically equivalent to a reducible \text{Fuchsian system}.
\newblock {\em Publication de L'Institut de Recherche Math\'ematique
  Avanc\'ee}, 2000.

\bibitem{malek3}
S.~Malek.
\newblock On the reducibility of the \text{Schlesinger} equations.
\newblock {\em J. Dynam. Control Systems}, \textbf{8}:no. 4, 505--527, 2002.

\bibitem{Mal}
B.~Malgrange.
\newblock {\em Sur les d\'eformations isomonodromiques I. Singularit\'es
  r\'eguli\`eres}, volume \textbf{37} of {\em Mathematics and Physics. Progr.
  Math.}
\newblock \text{Birkh\"auser} Boston, 1983.

\bibitem{Man1}
Yu.I. Manin.
\newblock {\em \text{Frobenius manifolds, quantum cohomology and moduli
  spaces}}, volume \textbf{47}.
\newblock American Mathematical Society, 1999.

\bibitem{mars-wein}
J.~Marsden and A.~Weinstein.
\newblock Reduction of symplectic manifolds with symmetry.
\newblock {\em Rep. Mathematical Phys.}, \textbf{5}:no. 1, 121--130, 1974.

\bibitem{M}
M.~Mazzocco.
\newblock Picard and \text{C}hazy solutions to the \text{PVI} equation.
\newblock {\em Math. Ann.}, \textbf{321}:no. 1, 131--169, 2001.

\bibitem{M2}
M.~Mazzocco.
\newblock Rational solutions of the \text{Painlev\'e VI} equation.
\newblock {\em J. Phys. A: Math. Gen.}, \textbf{34}:2281--2294, 2001.

\bibitem{M1}
M.~Mazzocco.
\newblock The geometry of the classical solutions of the \text{Garnier}
  systems.
\newblock {\em Int. Math. Res. Not.}, \textbf{2002}:no.12, 613--646, 2002.

\bibitem{Mi}
T.~Miwa.
\newblock Painlev\'e property of monodromy presereving equations and the
  analyticity of $\tau$-functions.
\newblock {\em Publ. RIMS}, \textbf{17}:703--721, 1981.

\bibitem{oht}
M.~Ohtsuki.
\newblock On the number of apparent singularities of a linear differential
  equation.
\newblock {\em Tokyo J. Math.}, \textbf{5}:23--29, 1982.

\bibitem{Ok1}
K.~Okamoto.
\newblock \text{Isomonodromic deformations, Painlev\'e} equations, and the
  \text{Garnier} system.
\newblock {\em J. Fac. Sci. Univ. Tokyo, Sect. 1A, Math.},
  \textbf{33}:576--618, 1986.

\bibitem{OK3}
K.~Okamoto.
\newblock Painlev\'e equations and \text{D}ynkin diagrams.
\newblock {\em Painlev\'e Transcendents}, pages 299--313, 1992.

\bibitem{OK}
K.~Okamoto and H.~Kimura.
\newblock \text{On Particular} solutions of the \text{Garnier} \text{system}
  and the \text{hypergeometric} \text{functions} of \text{several variables}.
\newblock {\em Quart. J. Math. Oxford}, \textbf{37}:61--80, 1986.

\bibitem{Sch}
L.~Schlesinger.
\newblock \text{Ueber eine Klasse von Differentsial System Beliebliger}
  \text{Ordnung} mit \text{Festen Kritischer Punkten}.
\newblock {\em J. fur Math.}, \textbf{141}:96--145, 1912.

\bibitem{sib}
Y.~Sibuya.
\newblock {\em Linear Differential Equations in the Complex Domain: Problems of
  Analytic Continuation}, volume \textbf{82}.
\newblock AMS TMM, 1990.

\bibitem{Um1}
H.~Umemura.
\newblock \text{Irreducebility of the first differential equation of
  Painlev\'e}.
\newblock {\em Nagoya Math. J.}, \textbf{117}:231--252, 1990.

\bibitem{singer}
Marius van~der Put and Michael~F. Singer.
\newblock {\em Galois theory of linear differential equations}, volume 328 of
  {\em Grundlehren der Mathematischen Wissenschaften [Fundamental Principles of
  Mathematical Sciences]}.
\newblock Springer-Verlag, Berlin, 2003.

\end{thebibliography}
\bibliographystyle{plain}

\end{document}